\documentclass[12pt]{article}
\usepackage{a4}
\usepackage{amsfonts}
\usepackage{amssymb}

\begin{document}

\title{\textbf{On the smallest poles of Igusa's p-adic zeta functions}}
\author{Dirk Segers}

\date{April 6, 2005}
\maketitle

\begin{abstract}
Let $K$ be a $p$-adic field. We explore Igusa's $p$-adic zeta
function, which is  associated to a  $K$-analytic function on an
open and compact subset of $K^n$. First we deduce a formula for an
important coefficient in the Laurent series of this meromorphic
function at a candidate pole. Afterwards we use this formula to
determine all values less than $-1/2$ for $n=2$ and less than $-1$
for $n=3$ which occur as the real part of a pole.
\end{abstract}

\section{Introduction}
\noindent \textbf{(1.1)} Let $K$ be a $p$-adic field, i.e., an
extension of $\mathbb{Q}_p$ of finite degree. Let $R$ be the
valuation ring of $K$, $P$ the maximal ideal of $R$ and $q$ the
cardinality of the residue field $R/P$. For $z \in K$, let
$\mathrm{ord} \, z \in \mathbb{Z} \cup \{+\infty\}$ denote the
valuation of $z$ and $|z|=q^{-\mathrm{ord} \, z}$ the absolute
value of $z$.

\vspace{0,5cm}

\noindent \textbf{(1.2)} Let $f$ be a $K$-analytic function on an
open and compact subset $X$ of $K^n$ and put $x=(x_1,\ldots,x_n)$.
Igusa's $p$-adic zeta function of $f$ is defined by
\[ Z_f(s)= \int_X |f(x)|^s \, |dx| \]
for $s \in \mathbb{C}$, $\mbox{Re}(s) \geq 0$, where $|dx|$
denotes the Haar measure on $K^n$, so normalized that $R^n$ has
measure $1$. Igusa proved that it is a rational function of
$q^{-s}$, so that it extends to a meromorphic function $Z_f(s)$ on
$\mathbb{C}$ which is also called Igusa's $p$-adic zeta function
of $f$.

\vspace{0,5cm}

\noindent \textbf{(1.3)} This zeta function has an interesting
connection with number theory. Let $f$ be a $K$-analytic function
on $R^n$ defined by a power series over $R$ which is convergent on
the whole of $R^n$. Let $M_i$ be the number of solutions of $f(x)
\equiv 0 \mbox{ mod } P^i$ in $(R/P^i)^n$. All the $M_i$'s are
described by $Z_f(s)$ through the relation

\vspace{0,5cm} \footnoterule{\footnotesize{2000 \emph{Mathematics
Subject Classification}. 11D79 11S80 14B05 14E15}}
\newpage

\noindent $Z_f(s) = (1-q^s)P(q^{-s})+q^s,$ where the Poincar\'e
series $P(t)$ of $f$ is defined by
\[ P(t) = \sum_{i=0}^{\infty} M_i(q^{-n}t)^i. \]
Remark that $P(t)$ is a rational function of $t$ because $Z_f(s)$
is a rational function of $q^{-s}$.

\vspace{0,5cm}

\noindent \textbf{(1.4)} The poles of $Z_f(s)$ are an interesting
object of study because they are related to the monodromy
conjecture \cite[(2.3.2)]{Denefreport} and because they determine
the asymptotic behaviour of the $M_i$. The poles with largest real
part give the largest contribution to the $M_i$. In this paper we
are concerned with the smallest real part $l$ of a pole of
$Z_f(s)$. A non-trivial consequence of the fact that the $M_i$ are
integers is that $l$ is larger than or equal to $-n$. Our main
results are stated in the next paragraph and sharpen this bound by
using a completely different method. This better bound has number
theoretic consequences because the knowledge of $l$ gives us
interesting information about the $M_i$: there exists an $a \in
\mathbb{Z}$ such that $M_i$ is divisible by
$q^{\ulcorner(n+l)i-a\urcorner}$ for all $i$ (for which $(n+l)i-a
\geq 0$). This is proved in the appendix. Remark that $a$ is
independent of $i$ and that the number in the exponent is the
smallest integer larger than or equal to $(n+l)i-a$.

Let $F_n^K$ denote the set of all $K$-analytic functions defined
on an arbitrary open and compact subset of $K^n$. For $n \in
\mathbb{Z}_{>0}$, we define the set $\mathcal{P}_n^K$ by
\[
\mathcal{P}_n^K := \{ s_0 \mid \exists f \in F_n^K \,
: \, Z_{f}(s) \textsl{\mbox{ has a pole with real part }} s_0 \}.
\]
In this article, we will prove that $\mathcal{P}_2^K \cap
]-\infty,-1/2[ = \{-1/2-1/i \mid i \in \mathbb{Z}_{>1} \} =
\{-1,-5/6,-3/4,-7/10,\ldots \}$ and that $\mathcal{P}_3^K \cap
]-\infty,-1[ = \{-1-1/i \mid i \in \mathbb{Z}_{>1} \}$. In
general, we expect that $\mathcal{P}_n^K \cap ]-\infty,-(n-1)/2[ =
\{-(n-1)/2-1/i \mid i \in \mathbb{Z}_{>1} \}$.

\vspace{0,2cm} \noindent \textsl{Remark.} One can easily show that
$\mathcal{P}_n^K \cap ]-\infty,-n+1[ = \emptyset$ if $n \geq 2$.

\vspace{0,5cm}

\noindent \textbf{(1.5)} Let $f \in K[x_1,x_2]$. Consider $f$ as a
polynomial over $K^{\mathrm{alg \, cl}}$. Suppose that the minimal
embedded resolution $g$ of $f^{-1}\{0\} \subset (K^{\mathrm{alg \,
cl}})^2$ is defined over $K$, i.e., all irreducible components of
$g^{-1}(f^{-1}\{0\})$ over $K^{\mathrm{alg \, cl}}$ and all points
in the intersection of two such components are defined over $K$.
Then it is generally known that an exceptional curve which is
intersected once or twice does not contribute to the residues of
its candidate poles with candidate order $1$. Because
$K^{\mathrm{alg \, cl}} \cong \mathbb{C}$, we can use the
calculations in \cite{SegersVeys} to conclude that the real part
of a pole of $Z_f(s)$ is of the form $-1/2-1/i$, $i \in
\mathbb{Z}_{>1}$, if it is smaller than $-1/2$.

Let $f \in K[x_1,x_2,x_3]$. Consider $f$ again as a polynomial
over $K^{\mathrm{alg \, cl}} \cong \mathbb{C}$. Suppose that there
exists an embedded resolution $g$ of $f^{-1}\{0\} \subset
(K^{\mathrm{alg \, cl}})^3 \simeq \mathbb{C}^3$ for which the
induced embedded resolution of the germ at each point $P$ of
$\mathbb{C}^3$ satisfies the conditions in
\cite[(3.1.1)]{SegersVeys}, which is defined over $K$ and which
has good reduction modulo $P$ (see \cite[section 2]{Denefdegree}).
Then the vanishing results in \cite{Veysconfigurations} and the
calculations in \cite{SegersVeys} imply that the real part of a
pole of $Z_f(s)$ is of the form $-1-1/i$, $i \in \mathbb{Z}_{>1}$,
if it is smaller than $-1$.

Consequently, starting from $\cite{SegersVeys}$, it is rather easy
to deal with polynomials which allow an appropriate embedded
resolution. However it is very difficult to verify the existence
of such an embedded resolution for a concrete function $f$. In a
lot of cases there does not exist an embedded resolution which is
defined over $K$, and if it exists, the condition of good
reduction modulo $P$ is very hard to check. This gives us a strong
motivation to study the general case. In this article there are no
constraints on $f$: we will not require that $g$ is defined over
$K$ and that $g$ has good reduction modulo $P$.

\vspace{0,5cm}

\noindent \textbf{Acknowledgements.} I want to thank Willem Veys
for advising me and for his useful remarks.

\section{The tool for our vanishing results}
\noindent \textbf{(2.1)} Let  $K$ be a $p$-adic field. Let $Y$ be
a $n$-dimensional $K$-analytic manifold, $\omega$ a $K$-analytic
differential $n$-form on $Y$ and $h$ a $K$-analytic function on
$Y$. We say that a chart $(V,y=(y_1,\ldots,y_n))$ on $Y$ is a good
chart for $(h,\omega)$ if $h= \varepsilon \prod_{i=1}^{k}
y_i^{N_i}$ on $V$ and $\omega = \eta \prod_{i=1}^{k} y_i^{\nu_i-1}
dy$ on $V$ for $k \in \{0,\ldots,n\}$, $N_i \in \mathbb{Z}_{> 0}$,
$\nu_i \in \mathbb{Z}_{>0}$ and non-vanishing $K$-analytic
functions $\varepsilon, \eta$ on $V$. We say that $(h,\omega)$ has
normal crossings at a point $P \in Y$ if there exists a good chart
for $(h,\omega)$ around $P$. So when we say normal crossings, we
mean normal crossings over $K$.

Let $x=(x_1,\ldots,x_n)$ be the coordinates of $K^n$. Let $f$ be a
$K$-analytic function on an open and compact subset $X$ of $K^n$.
Suppose that $f$ does not vanish on an open subset of $X$. An
embedded resolution of $(f,dx)$ consists of a $K$-analytic
$n$-dimensional manifold $Y$ and a proper $K$-analytic map $g:Y
\rightarrow X$ such that the restriction $Y \setminus
g^{-1}(f^{-1}\{0\}) \rightarrow X \setminus f^{-1}\{0\}$ is a
$K$-bianalytic map and such that $(f \circ g, g^*dx)$ has normal
crossings at every $P \in Y$. We can write $g^{-1}(f^{-1}\{0\})$
as a finite union of closed submanifolds $E_i$, $i \in T$, of
codimension one for which there exists a pair of positive integers
$(N_i,\nu_i)$, called the numerical data of $E_i$, such that the
following condition holds for every point $P \in Y$. If
$E_1,\ldots,E_k$ are all the $E_i$ that contain $P$, there exists
a chart $(V,y=(y_1,\ldots,y_n))$ around $P$ with $y_i$, $1 \leq i
\leq k$, an equation of $E_i$ on $V$ such that \vspace{-0,3cm}
\[
f \circ g = \varepsilon \prod_{i=1}^{k} y_i^{N_i} \qquad
\mbox{and} \qquad g^*dx= \eta \prod_{i=1}^{k} y_i^{\nu_i-1} dy
\vspace{-0.3cm}
\]
on $V$ for non-vanishing $K$-analytic functions $\varepsilon$ and
$\eta$ on $V$. By Hironaka's theorem \cite{Hironaka}, there always
exists an embedded resolution which is a composition of
blowing-ups along $K$-analytic closed submanifolds which are
contained in the zero locus of the pullback of $f$.

Let $g:Y \rightarrow X$ be a $K$-analytic map which is a
composition of blowing-ups along $K$-analytic closed submanifolds
which are contained in the zero locus of the pullback of $f$. If
$y=(y_1,\ldots,y_n)$ is a system of local parameters at $P \in Y$
such that $f \circ g = \varepsilon \prod_{i=1}^{k} y_i^{N_i}$ for
$k \in \{0,\ldots,n\}$, $N_i \in \mathbb{Z}_{>0}$ and
$\varepsilon$ a unit in the local ring at $P$, then $g^*dx = \eta
\prod_{i=1}^{k} y_i^{\nu_i-1}dy$ for $\nu_i \in \mathbb{Z}_{>0}$
and $\eta$ a unit in the local ring at $P$. Consequently we will
talk in this context about an embedded resolution of $f$, about
normal crossings of $f \circ g$ at $P$, and about a good chart for
$f \circ g$.

\vspace{0,5cm}

\noindent \textbf{(2.2)} Fix a uniformizing parameter $\pi$ for
$R$. For $z \in K$ let $\mbox{ac} \, z :=z \pi^{-\mathrm{ord} \,
z}$ be the angular component of $z$. Let $\chi$ be a character of
$R^{\times}$, i.e., a homomorphism $\chi: R^{\times} \rightarrow
\mathbb{C}^{\times}$ with finite image. Igusa's $p$-adic zeta
function of $f$ and $\chi$ is defined by
\[
Z_{f,\chi}(s)=\int_X \chi(\mathrm{ac} \, f(x)) |f(x)|^s \, |dx|
\]
for $s \in \mathbb{C}$, $\mathrm{Re}(s) \geq 0$. Note that
$Z_{f,\chi}(s)=Z_f(s)$ if $\chi$ is the trivial character, which
is denoted by 1.

Let $g:Y \rightarrow X$ be an embedded resolution of $(f,dx)$. We
study Igusa's $p$-adic zeta function $Z_{f,\chi}(s)$ by
calculating the integral on the resolution $Y$. Because
$|\varepsilon|$, $|\eta|$ and $\chi(\mathrm{ac} \, \varepsilon)$
are locally constant functions on each chart and because $Y$ is a
compact $K$-analytic manifold, we can choose a finite set $J$ of
good charts $(V,y)$ for $(f \circ g, g^*dx)$ such that $|
\varepsilon |$, $| \eta |$ and $\chi(\mathrm{ac} \, \varepsilon)$
are constant on each chart, the $V$'s form a partition of $Y$ and
for each chart $(V,y)$ we have $y(V)=P^j:=P^{j_1} \times \cdots
\times P^{j_n}$ for some $j=(j_1,\ldots,j_n) \in (\mathbb{Z}_{\geq
0})^n$ depending on $(V,y)$. Remark that we may even require that
$j_1= \cdots =j_n$ and that this value does not depend on the
chart, but we will not do this. We obtain
\begin{eqnarray*}
Z_{f,\chi}(s) & = & \int_X \chi(\mathrm{ac} \, f(x)) |f(x)|^s \,
|dx| \\ & = & \sum_{(V,y) \in J} \int_V \chi(\mathrm{ac} \,
(f\circ g)(y)) |(f \circ g)(y)|^s \, |g^*dx|
\\ & = & \sum_{(V,y) \in J} \int_{P^j} \chi(\mathrm{ac} \,
\varepsilon \textstyle{\prod_{i=1}^{k}} y_i^{N_i}) |\varepsilon
\textstyle{\prod_{i=1}^{k}} y_i^{N_i}|^s \, |\eta
\textstyle{\prod_{i=1}^{k}} y_i^{\nu_i-1} dy| \\& = & \sum_{(V,y)
\in J} \chi(\mathrm{ac} \, \varepsilon) |\varepsilon|^s |\eta|
q^{-\sum_{i=k+1}^{n}j_i} \prod_{i=1}^{k} \int_{P^{j_i}}
\chi^{N_i}(\mathrm{ac} \, y_i) |y_i|^{N_is+\nu_i-1} \, |dy_i|.
\end{eqnarray*}
The remaining integrals are easy and well known: in the one
variable case, we have that
\[
\int_{P^j} \chi(\mathrm{ac} \, x) |x|^{\alpha-1} \, |dx| = \left\{
\begin{array}{ll}
\frac{q-1}{q} \frac{q^{-j\alpha}}{1-q^{-\alpha}} & \mathrm{if} \,
\chi=1  \\ 0 & \mathrm{if} \, \chi \not= 1
\end{array} \right.
\]
for a complex number $\alpha$ with $\mathrm{Re}(\alpha) > 0$ and
that the integral is not defined if $\mathrm{Re}(\alpha) \leq 0$.
This shows that $Z_{f,\chi}(s)$ is a rational function of
$q^{-s}$, so that it extends to a meromorphic function
$Z_{f,\chi}(s)$ on $\mathbb{C}$. Moreover, these calculations
imply that the integral $\int_X \chi(\mathrm{ac} \, f(x)) |f(x)|^s
|dx|$ is defined if and only if $\mbox{Re}(s) > \max \{ -\nu_i/N_i
\mid i \in T \}$. We obtain also from this calculation that every
pole of $Z_{f,\chi}(s)$ is of the form
\[
-\frac{\nu_i}{N_i} + \frac{2k \pi \sqrt{-1}}{N_i \log q},
\]
with $k \in \mathbb{Z}$ and $i \in T$ such that $\chi^{N_i}=1$.
These values are called the candidate poles of $Z_{f,\chi}(s)$.
The candidate poles of $Z_{f,\chi}(s)$ associated to $E_i$, with
$i \in T$ such that $\chi^{N_i}=1$, are the values $-\nu_i/N_i +
(2k \pi \sqrt{-1})/(N_i \log q)$, $k \in \mathbb{Z}$. Obviously,
we do not associate candidate poles to $E_i$ if $\chi^{N_i} \not=
1$.

Let $s_0$ be a candidate pole of $Z_{f,\chi}(s)$. Because the
poles of $1/(1-q^{-N_is- \nu_i})$ have order one, we define the
expected order $m=m(s_0)$ of $s_0$ as the highest number of
$E_i$'s with candidate pole $s_0$ and with non-empty intersection.
The order of $s_0$ is of course less than or equal to $m$. It is
less than $m$ if and only if $b_{-m}$, which is defined by the
Laurent series
\[ \frac{b_{-m}}{(s-s_0)^m} + \frac{b_{-m+1}}{(s-s_0)^{m-1}}+ \cdots + b_0 + b_1(s-s_0) + \cdots \]
of $Z_{f,\chi}(s)$ at $s_0$, is equal to zero. Remark that a
candidate pole of expected order one is a pole if and only if
$b_{-1} \not= 0$.

Everything we have done up till now is well known. More details
can be found for example in \cite{Igusabook}.

\vspace{0,5cm}

\noindent \textbf{(2.3)} Let $X$ be an open and compact subset of
$K^n$. Let $\xi, f, f_1, \ldots, f_l$ be $K$-analytic functions on
$X$. Let $a_i,b_i$, $1 \leq i \leq l$, be non-negative integers.
We associate to these data the zeta function
\[
Z(s_1,\ldots,s_l) = \int_X \chi(\mathrm{ac} \, f) |\xi|
|f_1|^{a_1s_1+b_1} \ldots |f_l|^{a_ls_l+b_l} |dx|,
\]
which is defined on a set $U$ that contains all points
$(s_1,\ldots,s_l) \in \mathbb{C}^l$ with $\mbox{Re}(s_i) \geq
-b_i/a_i$ if $f_i$ vanishes on $X$ and $s_i$ arbitrary if $f_i$
does not vanish on $X$. Loeser \cite{Loeser2} already studied this
zeta function. By looking at an embedded resolution of $\xi f f_1
\ldots f_l$, one proves in a way that is analogous to the argument
in (2.2) that $Z(s_1,\ldots,s_l)$ is a rational function of
$q^{-s_1}, \ldots, q^{-s_l}$. Consequently, it extends to a
meromorphic function on $\mathbb{C}^l$, which we also denote by
$Z(s_1,\ldots,s_l)$. As before, we can also obtain an explicit
description of $U$, which turns out to be an open subset of
$\mathbb{C}^l$.

The meromorphic continuation of a function $h$ will be denoted by
$[ h ]^{\mathrm{mc}}$ and the evaluation of this meromorphic
continuation at the point $s=s_0$ of the domain will be denoted by
$[ h ]^{\mathrm{mc}}_{s=s_0}$.

In our study of Igusa's $p$-adic zeta function, we will have to
deal with expressions of the form
\[
\left[ \int_X \chi(\mathrm{ac} \, f) |\xi| |f_1|^{a_1s+b_1} \ldots
|f_l|^{a_ls+b_l} |dx| \right]^{\mathrm{mc}}_{s=s_0}.
\]
The zeta function in more complex variables can be used to modify
this expression. If $U \cap \{ (s_1,\ldots,s_l) \in \mathbb{C}^l
\mid s_1 = s_0 \} \not= \emptyset$, then
\begin{eqnarray*}
\lefteqn{ \left[ \int_X \chi(\mathrm{ac} \, f) |\xi|
|f_1|^{a_1s+b_1} |f_2|^{a_2s+b_2} \ldots |f_l|^{a_ls+b_l} |dx|
\right]^{\mathrm{mc}}_{s=s_0} } \\ & = & \left[ \int_X
\chi(\mathrm{ac} \, f) |\xi| |f_1|^{a_1s_1+b_1} |f_2|^{a_2s_2+b_2}
\ldots |f_l|^{a_ls_l+b_l} |dx|
\right]^{\mathrm{mc}}_{s_1=\cdots=s_l=s_0}
\\ & = & \left[ \int_X \chi(\mathrm{ac} \, f) |\xi| |f_1|^{a_1s_0+b_1} |f_2|^{a_2s+b_2}
\ldots |f_l|^{a_ls+b_l} |dx| \right]^{\mathrm{mc}}_{s=s_0}.
\end{eqnarray*}
We explain the first equality. The composition of the map $A :
\mathbb{C} \rightarrow \mathbb{C}^{l} : s \mapsto (s,\ldots,s)$
with the meromorphic function $Z(s_1,\ldots,s_l)$ on
$\mathbb{C}^{l}$ is a meromorphic function on $\mathbb{C}$ which
is equal to the meromorphic function
\[
\left[ \int_X \chi(\mathrm{ac} \, f) |\xi| |f_1|^{a_1s+b_1} \ldots
|f_l|^{a_ls+b_l} |dx| \right]^{\mathrm{mc}}
\]
because they agree on an open subset of $\mathbb{C}$.
Consequently, the first equality is nothing more than $(Z \circ
A)(s_0)=Z(A(s_0))$. For the second equality, we have to use the
map $B : \mathbb{C} \rightarrow \mathbb{C}^{l} : s \mapsto
(s_0,s,\ldots,s)$.

\vspace{0,5cm}

\noindent \textbf{(2.4)} Let $f$ be a $K$-analytic function on an
open and compact subset $X$ of $K^n$ and let $g:Y \rightarrow X$
be an embedded resolution of $(f,dx)$. Denote $E_I= \cap_{i \in I}
E_i$ for $I \subset T$. Let $\chi$ be a character of $R^{\times}$.
Let $s_0$ be a candidate pole of $Z_{f,\chi}(s)$ and let $m$ be
its expected order. Let $E_I$, $I \in S$, be all the non-empty
intersections of $m$ varieties $E_i$, $i \in T$, with candidate
pole $s_0$ (and thus also with $\chi^{N_i}=1$). Fix $I \in S$ and
suppose for the ease of notation that $I=\{1,\ldots,m\}$. Let
$W_1$ and $W_2$ be open and compact subsets of $Y$ which satisfy
$E_I \cap W_1 = E_I \cap W_2 \not= \emptyset$ and which do not
meet any $E_K$, $K \in S \setminus \{I\}$. Then the contribution
of $W_1$ to $b_{-m}$ and the contribution of $W_2$ to $b_{-m}$ are
the same because they are both equal to the contribution of $W_1
\cap W_2$ to $b_{-m}$. Consequently we can speak of the
contribution of $E_I \cap W_1 = E_I \cap W_2$ to $b_{-m}$. In
particular, the contribution of $E_I$ to $b_{-m}$ is well defined.

Consider a set $J$ of disjoint compact charts $(V,y)$ that
intersect $E_I$, that cover $E_I$ and that are disjoint with all
$E_K$, $K \in S \setminus \{I\}$. This set $J$ is necessarily
finite and the contribution of $E_I$ to $b_{-m}$ is the sum over
$J$ of the contributions
\[
\lim_{s \rightarrow s_0} (s-s_0)^m \left[ \int_V \chi(\mathrm{ac}
\, (f \circ g)(y)) |(f \circ g)(y)|^s |g^* dx|
\right]^{\mathrm{mc}}
\]
of $V$ to $b_{-m}$.

We introduce some notation. Let $(V,y)$ be a chart. We have that
$\overline{y}=(y_{m+1},\ldots,y_n)$ determines a chart on the
closed submanifold $\overline{V}$ defined by $y_1= \cdots =
y_m=0$. Denote $dy_{m+1} \wedge \cdots \wedge dy_n$ by
$d\overline{y}$. It is a volume form on $\overline{V}$. If
$j=(j_1, \ldots ,j_n) \in (\mathbb{Z}_{\geq 0})^n$, then we denote
$P^{j_{m+1}} \times \cdots \times P^{j_n}$ by $\overline{P^j}$.

Suppose that $(V,y)$ is a chart such that $E_1,\ldots,E_m$ have
equations $y_1=0, \ldots, y_m=0$ respectively, and such that
\[
f \circ g = \alpha \prod_{i=1}^{m} y_i^{N_i} \qquad \mbox{and}
\qquad g^*dx = \beta \prod_{i=1}^{m} y_i^{\nu_i-1} dy
\]
on $V$, for $K$-analytic functions $\alpha$ and $\beta$ on $V$
with $|\alpha|$, $|\beta|$ and $\chi(\mathrm{ac} \, \alpha)$
\textsl{independent of} $y_1,\ldots,y_m$. Remark that a good chart
$(V,y)$ for $(f \circ g,g^*dx)$ in which $|\varepsilon|$, $|\eta|$
and $\chi(\mathrm{ac} \, \varepsilon)$ are constant satisfies this
condition for $\alpha = \varepsilon \prod_{i=m+1}^{k} y_i^{N_i}$
and $\beta= \eta \prod_{i=m+1}^{k} y_i^{\nu_i-1}$. Remark also
that $\overline{V} = V \cap E_I$. Suppose also that $y(V)$ is of
the form $P^j$ with $j=(j_1,\ldots,j_n) \in (\mathbb{Z}_{\geq
0})^n$. Then
\begin{eqnarray*}
\lefteqn{ \lim_{s \rightarrow s_0} (s-s_0)^m \left[ \int_V
\chi(\mathrm{ac} \, (f \circ g)(y)) |(f \circ g)(y)|^s |g^*dx|
\right]^{\mathrm{mc}} } \\ & = & \lim_{s \rightarrow s_0}
(s-s_0)^m \left[ \int_{P^j} \chi(\mathrm{ac} \, \alpha) |\alpha|^s
|\beta| \prod_{i=1}^{m} \chi^{N_i}(\mathrm{ac} \, y_i)
|y_i|^{N_is+ \nu_i -1} |dy| \right]^{\mathrm{mc}}
\\ & = & \left( \prod_{i=1}^{m} \lim_{s
\rightarrow s_0} (s-s_0) \left[ \int_{P^{j_i}} |y_i|^{N_is+ \nu_i
-1} |dy_i| \right]^{\mathrm{mc}} \right) \left[
\int_{\overline{P^j}} \chi(\mathrm{ac} \, \alpha) |\alpha|^s
|\beta| |d\overline{y}| \right]^{\mathrm{mc}}_{s=s_0}
\\ & = & \left(\prod_{i=1}^{m} \frac{q-1}{qN_i \log q}
\right) \left[ \int_{\overline{V}} \chi(\mathrm{ac} \, \alpha)
|\alpha|^s |\beta| |d\overline{y}| \right]^{\mathrm{mc}}_{s=s_0}
\end{eqnarray*}

We have shown that the last expression is the contribution of $V$
to $b_{-m}$. Consequently, the only aspect of the chart $(V,y)$
that it depends on is $\overline{V}$. In the next section we will
see that this is still the case if $|\alpha|$, $|\beta|$ and
$\chi(\mathrm{ac} \, \alpha)$ depend on $y_1,\ldots,y_m$ and if we
are not in an embedded resolution.

\vspace{0,5cm}

\noindent \textbf{(2.5)} Suppose that $g:Y=Y_t \rightarrow X=Y_0$
is a composition $g_1 \circ \cdots \circ g_t$ of blowing-ups
$g_i:Y_i \rightarrow Y_{i-1}$. Suppose that each $g_i$ is a
blowing-up along a $K$-analytic closed submanifold of codimension
bigger than one which has only normal crossings with the union of
the exceptional varieties of $g_1 \circ \cdots \circ g_{i-1}$. Let
$I=\{1,\ldots,m\} \in S$ as in (2.4). Let $r \in \{ 0,\ldots,t
\}$. Suppose that $E_I$ already exists in $Y_r$ and that the
$E_i$, $i \in I$, intersect transversally in $Y_r$. Remark that
the last condition is satisfied if all the $E_i$, $i \in I$, are
exceptional. We will write $E_I \subset Y_r$ if we want to stress
that we consider $E_I$ as a subset of $Y_r$.

We call a chart $(V,y)$ a good chart for $E_I \subset Y_r$ if
$(V,y)$ is a chart on $Y_r$ such that $V$ intersects $E_I$ and
such that $y_1=0,\ldots,y_m=0$ are the equations of respectively
$E_1,\ldots,E_m$ on $V$.

Let $(V,y)$ be a good chart for $E_I \subset Y_r$. Then we have
\[
f \circ g_1 \circ \cdots \circ g_r = \alpha \prod_{i=1}^{m}
y_i^{N_i} \quad \mbox{and} \quad (g_1 \circ \cdots \circ g_r)^* dx
= \beta \prod_{i=1}^{m} y_i^{\nu_i-1} dy
\]
on $V$, for $K$-analytic functions $\alpha$ and $\beta$ on $V$.

\vspace{0,5cm}

We will now prove that the only aspect of the chart $(V,y)$ that
\begin{eqnarray}
\left[ \int_{\overline{V}} \chi(\mathrm{ac} \, \alpha) |\alpha|^s
|\beta| |d\overline{y}| \right]^{\mathrm{mc}}_{s=s_0}
\end{eqnarray}
depends on is $\overline{V}$.

Let $(W,z)$ be another chart on $Y_r$ such that
$\overline{V}=\overline{W}$ and such that $z_1=0,\ldots,z_m=0$ are
the equations of respectively $E_1,\ldots,E_m$ on $W$. We may
suppose that $V=W$ because we can restrict them both to $V \cap
W$. For every $i \in \{ 1, \ldots ,m \}$ there exists a
non-vanishing $K$-analytic function $f_i$ on $V$ such that
$y_i=f_iz_i$ because $y_i$ and $z_i$ are equations of the same
$E_i$. Thus
\[
f \circ g_1 \circ \cdots \circ g_r = \alpha \left( \prod_{i=1}^m
f_i^{N_i} \right) \prod_{i=1}^m z_i^{N_i} \] and
\[
(g_1 \circ \cdots \circ g_r)^* dx = \beta \left( \prod_{i=1}^m
f_i^{\nu_i-1} \right) \det \left( \frac{\partial y}{\partial z}
\right) \prod_{i=1}^m z_i^{\nu_i-1} dz.
\]
We have to prove that (1) is equal to
\begin{eqnarray}
\left[ \int_{\overline{V}} \chi(\mathrm{ac} \, \alpha) |\alpha|^s
|\beta| \left( \prod_{i=1}^{m} \chi^{N_i}(\mathrm{ac} \, f_i)
|f_i|^{N_is+ \nu_i -1} \right) \left| \det \left( \frac{\partial
y}{\partial z} \right) \right| |d\overline{z}|
\right]^{\mathrm{mc}}_{s=s_0}.
\end{eqnarray}
In (1) we have that $|d\overline{y}| = |\det (\partial
\overline{y}/ \partial \overline{z})| |d\overline{z}|$. Recall
that $\chi^{N_i}=1$ for $i \in \{ 1,\ldots,m \}$. Because $f_i$,
$i \in I$, is a non-vanishing function, we may replace each
$|f_i|^{N_i s + \nu_i -1}$ in (2) by $|f_i|^{N_i s_0 + \nu_i -1}$
according to (2.3), and this is equal to $|f_i|^{-1}$ because $N_i
s_0 + \nu_i =0$ for $i \in I$. Consequently, we have to prove that
\begin{eqnarray*}
\lefteqn{ \left[ \int_{\overline{V}} \chi(\mathrm{ac} \, \alpha)
|\alpha|^s |\beta| \left| \det \left( \frac{\partial
\overline{y}}{\partial \overline{z}} \right) \right|
|d\overline{z}| \right]^{\mathrm{mc}}_{s=s_0} } \\ & = & \left[
\int_{\overline{V}} \chi(\mathrm{ac} \, \alpha) |\alpha|^s |\beta|
\prod_{i=1}^{m} |f_i|^{-1} \left| \det \left( \frac{\partial
y}{\partial z} \right) \right| |d\overline{z}|
\right]^{\mathrm{mc}}_{s=s_0}.
\end{eqnarray*}
Because $(\partial y / \partial z)$ is equal to
\[
\left( \begin{array}{cccccc} f_1+z_1 \frac{\partial f_1}{\partial
z_1} & \cdots & z_1 \frac{\partial f_1}{\partial z_m} & z_1
\frac{\partial f_1}{\partial z_{m+1}} & \cdots & z_1
\frac{\partial f_1}{\partial z_n} \\ \vdots & \ddots &  \vdots &
\vdots & \ddots & \vdots \\ z_m \frac{\partial f_m}{\partial z_1}
& \cdots & f_m + z_m \frac{\partial f_m}{\partial z_m} & z_m
\frac{\partial f_m}{\partial z_{m+1}} & \cdots & z_m
\frac{\partial f_m}{\partial z_n} \\ \frac{\partial
y_{m+1}}{\partial z_1} & \cdots & \frac{\partial y_{m+1}}{\partial
z_m} & \frac{\partial y_{m+1}}{\partial z_{m+1}} & \cdots &
\frac{\partial y_{m+1}}{\partial z_n}
\\ \vdots & \ddots & \vdots & \vdots & \ddots & \vdots \\ \frac{\partial y_n}{\partial z_1} &
\cdots & \frac{\partial y_n}{\partial z_m} & \frac{\partial
y_n}{\partial z_{m+1}} & \cdots & \frac{\partial y_n}{\partial
z_n}
\end{array} \right) ,
\]
we obtain that
\[
\left[ \det \left( \frac{\partial y}{\partial z} \right)
\right]_{z_1= \cdots =z_m=0} = \left[ \left( \prod_{i=1}^{m} f_i
\right) \det \left( \frac{\partial \overline{y}}{\partial
\overline{z}} \right) \right]_{z_1= \cdots =z_m=0}.
\]
Consequently we have proved our statement.

\vspace{0,5cm}

\noindent \textbf{(2.6)} In order to formulate the next important
proposition, we recall the setting.

\vspace{0,2cm}

\noindent \textsl{Data and notations.}  Let $f$ be a $K$-analytic
function on an open and compact subset $X$ of $K^n$. Let $\chi$ be
a character of $R^{\times}$. Let $g:Y=Y_t \rightarrow X=Y_0$ be an
embedded resolution of $(f,dx)$ which is a composition $g_1 \circ
\cdots \circ g_t$ of blowing-ups $g_i:Y_i \rightarrow Y_{i-1}$.
Suppose that each $g_i$ is a blowing-up along a $K$-analytic
closed submanifold of codimension larger than one which has only
normal crossings with the union of the exceptional varieties of
$g_1 \circ \cdots \circ g_{i-1}$. Let $s_0$ be a candidate pole of
$Z_{f,\chi}(s)$ and let $m$ be its expected order. Let $b_{-m}$ be
the coefficient of $1/(s-s_0)^m$ in the Laurent series of
$Z_{f,\chi}(s)$ at $s_0$. For $I \subset T$, denote the
intersection in $Y$ of the $E_i \subset Y$, $i \in I$, by $E_I$.
Let $E_I$, $I \in S$, be all the non-empty intersections in $Y$ of
$m$ varieties $E_i$, $i \in T$, with candidate pole $s_0$ (and
thus also with $\chi^{N_i}=1$). Fix $I \in S$ and suppose for the
ease of notation that $I=\{1,\ldots,m\}$. Let $r \in \{ 0,\ldots,t
\}$. Suppose that $E_I$ already exists in $Y_r$ and that the
$E_i$, $i \in I$, intersect transversally in $Y_r$.

\vspace{0,2cm}

\noindent \textbf{Definition.} The contribution of an open and
compact subset $U$ of $E_I \subset Y_r$ to $b_{-m}$ is the
contribution of the strict transform of $U$ in $Y$ to $b_{-m}$.

\vspace{0,2cm}

\noindent \textsl{Remark.} The contribution of $E_I$ to $b_{-m}$
is not necessarily equal to the contribution of a `very small'
neighbourhood of $E_I \subset Y_r$ to $b_{-m}$, because it can
happen that an $E_K$, $K \in S \setminus \{I\}$, lies above $E_I
\subset Y_r$.

\vspace{0,2cm}

\noindent \textsl{Data and notations.} Let $(V,y)$ be a good
compact chart for $E_I \subset Y_r$, i.e., $(V,y)$ is a good chart
on $Y_r$ such that $V$ is compact. Write
\[
f \circ g_1 \circ \cdots \circ g_r = \alpha \prod_{i=1}^{m}
y_i^{N_i} \quad \mbox{and} \quad (g_1 \circ \cdots \circ g_r)^* dx
= \beta \prod_{i=1}^{m} y_i^{\nu_i-1} dy
\]
on $V$, for $K$-analytic functions $\alpha$ and $\beta$ on $V$. We
have that $\overline{y}=(y_{m+1},\ldots,y_n)$ determines
coordinates on the closed submanifold $\overline{V}=V \cap E_I$
which is defined by $y_1=\cdots=y_m=0$. Consider the volume form
$d\overline{y}=dy_{m+1} \wedge \cdots \wedge dy_n$ on
$\overline{V}$.

\vspace{0,2cm}

\noindent \textbf{Proposition.} \textsl{The contribution of
$\overline{V}$ to $b_{-m}$ is equal to}
\[
\left(\prod_{i=1}^{m} \frac{q-1}{qN_i \log q} \right) \left[
\int_{\overline{V}} \chi(\mathrm{ac} \, \alpha) |\alpha|^s |\beta|
|d\overline{y}| \right]^{\mathrm{mc}}_{s=s_0}.
\]

\vspace{0,2cm}

\noindent \emph{Notation.} For the ease of notation, we denote
$\kappa = \prod_{i=1}^{m} (q-1)/(q N_i \log q)$.

\vspace{0,2cm}

\noindent \emph{Proof.} We proved it already in the case $r=t$.
Suppose from now on that $r<t$.

Let $F$ be the closed subset of $E_I \subset Y_t$ that contains
the points at which $\pi = g_t \circ \cdots \circ g_{r+1} : Y_t
\rightarrow Y_r$ is not $K$-bianalytic. Let $(B_n)_{n \in
\mathbb{Z}_{\geq 0}}$ be a decreasing sequence of open and compact
subsets of $E_I \subset Y_t$ such that $\cap_{n \in
\mathbb{Z}_{\geq 0}} B_n = F$. Note that $\pi(F)$ is a closed
subset of $Y_r$ and that $\cap_{n \in \mathbb{Z}_{\geq 0}}
\pi(B_n) = \pi(F)$. Note also that for every neighbourhood $O$ of
$F$ (respectively $\pi(F)$), there exists a positive integer $n$
such that $B_n \subset O$ (respectively $\pi(B_n) \subset O$).

Let $J$ be a set of good compact charts $(W,z)$ for $E_I \subset
Y_t$ such that the $\overline{W}$'s form a partition of
$\pi^{-1}(\overline{V})$. Let $s \in \mathbb{C}$ with
$\mathrm{Re}(s)>0$. Then
\begin{eqnarray*}
\int_{\overline{V}} \chi(\mathrm{ac} \, \alpha) |\alpha|^s |\beta|
|d\overline{y}| & = & \lim_{n \rightarrow \infty}
\int_{\overline{V} \setminus \pi(B_n)} \chi(\mathrm{ac} \, \alpha)
|\alpha|^s |\beta| |d\overline{y}| \\ & = & \lim_{n \rightarrow
\infty} \sum_{(W,z) \in J} \int_{\overline{W} \setminus B_n}
\chi(\mathrm{ac} \, \alpha_z) |\alpha_z|^s |\beta_z|
|d\overline{z}| \\ & = & \sum_{(W,z) \in J} \int_{\overline{W}}
\chi(\mathrm{ac} \, \alpha_z) |\alpha_z|^s |\beta_z|
|d\overline{z}|.
\end{eqnarray*}
The first equality holds because
\[
\lim_{n \rightarrow \infty} \int_{\overline{V} \cap \pi(B_n)}
\chi(\mathrm{ac} \, \alpha) |\alpha|^s |\beta| |d\overline{y}| =
0.
\]
Indeed, the measure of $\overline{V} \cap \pi(B_n)$ for
$|d\overline{y}|$ decreases to zero if $n \rightarrow \infty$ and
the real and the complex part of the integrand are bounded for
complex numbers $s$ satisfying $\mathrm{Re}(s)>0$. The last
equality is obtained by using the same argument. We have written
$\alpha_z$ and $\beta_z$ to stress that these functions depend on
the chart. For the second equality, we have to use that the
$\overline{W}$'s form a partition of $\pi^{-1}(\overline{V})$,
that $\pi : Y_t \rightarrow Y_r$ is a $K$-bianalytic map on a
neighbourhood of $\overline{W} \setminus B_n$ in $Y_t$ for every
$(W,z) \in J$ and that the contribution is independent of the
chosen coordinates, a fact we explained in (2.5).

Finally, we evaluate the meromorphic continuation in $s=s_0$ and
we obtain:
\[
\left[ \int_{\overline{V}} \chi(\mathrm{ac} \, \alpha) |\alpha|^s
|\beta| |d\overline{y}| \right]^{\mathrm{mc}}_{s=s_0} =
\sum_{(W,z) \in J} \left[ \int_{\overline{W}} \chi(\mathrm{ac} \,
\alpha_z) |\alpha_z|^s |\beta_z| |d\overline{z}|
\right]^{\mathrm{mc}}_{s=s_0}.
\]
The right hand side multiplied by $\kappa$ is the contribution of
$\overline{V}$ to $b_{-m}$ because the proposition is true in the
case $r=t$. Consequently, the proof is done. $\qquad \Box$

\vspace{0,2cm}

\noindent \textsl{Remark.}  One can also prove that the
proposition is true on $Y_r$ if it is true on $Y_{r+1}$ by using
coordinate transformations of the blowing-up $g_{r+1}$. This
alternative proof can be found in \cite[page 58]{Segers}.

\vspace{0,5cm}

\noindent \textbf{(2.7)} Let $T_t$ be the set of all $j \in T
\setminus I$ for which $E_j$ intersects $E_I$ in $Y_t$. Let $F_j$,
$j \in T_t$, be the intersection of $E_j$ and $E_I$ in $Y_t$. We
have that $F_j$ has codimension one in $E_I \subset Y_t$. The set
of all $j$, $j \in T_t$, for which $(g_{r+1} \circ \cdots \circ
g_t)(F_j)$ has also codimension one in $E_I \subset Y_r$ will be
denoted by $T_r$. For $j \in T_r$ we denote $(g_{r+1} \circ \cdots
\circ g_t)(F_j)$ also by $F_j$ and we put $\alpha_j=N_js_0+\nu_j$.

Let $(V,y)$ be a good chart for $E_I$ on $Y_r$ on which $F_j$, $j
\in T_r$, is given by $y_1= \cdots = y_m = y_{m+1}=0$. Write
\begin{eqnarray*}
\alpha(0,\ldots,0,y_{m+1},\ldots,y_n) & = & y_{m+1}^{N_{j,r}} h_1
\qquad \mbox{and} \\ \beta(0,\ldots,0,y_{m+1},\ldots,y_n) & = &
y_{m+1}^{\nu_{j,r}-1} h_2
\end{eqnarray*}
with $h_1$ and $h_2$ not divisible by $y_{m+1}$. Then we denote
$N_{j,r}s_0+\nu_{j,r}$ by $\alpha_{j,r}$.

We deduce now the relations that will be used later. In this
paragraph we suppose that $m=1$ and that $E_I=E_r$ is created by
the blowing-up at a point $P$ of $Y_{r-1}$. Suppose that there
exists a chart $(V,y)$ centred at $P$ on which $f \circ g_1 \circ
\cdots \circ g_{r-1}$ is given by a power series with lowest
degree part a homogeneous polynomial for which every irreducible
factor over $K^{\mathrm{alg \, cl}}$ is defined over $K$ and for
which the zero locus in $\mathbb{P}^{n-1}$ of every irreducible
factor (over $K^{\mathrm{alg \, cl}}$) contains a non-singular
point defined over $K$. Remark that these conditions are satisfied
if the lowest degree part is a product of linear factors defined
over $K$. Write $f \circ g_1 \circ \cdots \circ g_{r-1}=e \left(
\prod_{j \in T_r} f_j^{N_{j,r}} \right) + \theta$ and $(g_1 \circ
\cdots \circ g_{r-1})^*dx=\rho \left( \prod_{j \in T_r}
f_j^{\nu_{j,r}-1} \right) dy$, where $f_j$ is the equation of $F_j
\subset E_r$ in the homogeneous coordinates $(y_1:\cdots:y_n)$ on
$E_r \subset Y_r$, $e \in K^{\times}$, $\theta$ is a power series
with multiplicity larger than the degree of the homogeneous
polynomial $\prod_{j \in T_j} f_j^{N_{j,r}}$ and $\rho$ is a
$K$-analytic function which does not vanish at $P$. Because the
multiplicity of $f \circ g_1 \circ \cdots \circ g_{r-1}$ at $P$ is
equal to $N_r$, we obtain the first relation:
\[
\sum_{j \in T_r} (\deg F_j) N_{j,r} = N_r.
\]
\vspace{-0,2cm}
\begin{picture}(150,0)(0,0)
\put(115,10){(\textsc{Relation 1})}
\end{picture}
Our second relation will involve the $\alpha_{j,r}$, $j \in T_r$.
There will appear differential forms with rational exponents in
the calculations. One can make sense to this by considering them
as an element of a tensor power of the module of rational
differential forms (see \cite{Jacobs}), but we will not give
details here. Let $i \in \{1,\ldots,n\}$. We look at the chart
$(O,z=(z_1,\ldots,z_n))$ on $Y_r$ for which
$g_r(z_1,\ldots,z_n)=(z_1z_i,\ldots,z_{i-1}z_i,z_i,z_{i+1}z_i,\ldots,z_nz_i)$.
Then
\[
f \circ g_1 \circ \cdots \circ g_r = z_i^{N_r} \left( e \prod_{j
\in T_r} f_j(z_1,\ldots,z_{i-1},1,z_{i+1},\ldots,z_n)^{N_{j,r}} +
z_i \frac{\theta \circ g_r}{z_i^{N_r+1}} \right)
\]
and
\[
(g_1 \circ \cdots \circ g_r)^*dx = z_i^{\nu_r-1} (\rho \circ g_r)
\left( \prod_{j \in T_r}
f_j(z_1,\ldots,z_{i-1},1,z_{i+1},\ldots,z_n)^{\nu_{j,r}-1} \right)
dz.
\]
Consequently the Poincar\'e residue of $(f \circ g_1 \circ \cdots
\circ g_r)^{-\nu_r/N_r}(g_1 \circ \cdots g_r)^*dx$ on $E_I \subset
Y_r$ (see \cite{Jacobs}) is equal to
\[ e^{-\nu_r/N_r} \rho(P) \prod_{j \in T_r}
f_j(z_1,\ldots,z_{i-1},1,z_{i+1},\ldots,z_n)^{\alpha_{j,r}-1}d\overline{z},
\]
so that the canonical divisor of $E_r$ is $\sum_{j \in T_r}
(\alpha_{j,r}-1)F_j$. Because we know that the degree of the
canonical divisor on $E_r \cong \mathbb{P}^{n-1}$ is $-n$, we
obtain the second relation: \vspace{-0,4cm}
\begin{picture}(150,0)(0,0)
\put(115,-5){(\textsc{Relation 2})}
\end{picture}
\[
\sum_{j \in T_r} (\deg F_j)(\alpha_{j,r}-1)=-n.
\]
Remark that the condition on the lowest degree part of $f \circ
g_1 \circ \cdots \circ g_{r-1}$ has to be satisfied because
otherwise some terms on the left hand side are missing. We need
the two relations which we just derived in section 3. In the next
paragraph we will deduce that $\alpha_{j,r}=\alpha_j$ and that
$N_{j,r} \equiv N_j \pmod{N_r}$ so that we obtain
\begin{eqnarray*}
& & \sum_{j \in T_r} (\deg F_j) N_j \equiv 0 \pmod{N_r} \qquad
\mbox{and} \\ & & \sum_{j \in T_r} (\deg F_j) (\alpha_j -1) = -n.
\end{eqnarray*}
One can find these relations in a more general form in
\cite{Veysrelations}, \cite{Veyscongruences1} and
\cite{Veyscongruences2}.

We prove that $\alpha_{j,r}=\alpha_j$ for $j \in T_r$. Because
$g_{r+1} \circ \cdots \circ g_t$ is a composition of a finite
number of blowing-ups, it is enough to prove that
$\alpha_{j,r}=\alpha_{j,r+1}$. If the centre of $g_{r+1}$ does not
contain $F_j$, then $N_{j,r}=N_{j,r+1}$ and
$\nu_{j,r}=\nu_{j,r+1}$ so that we are done. If the centre of
$g_{r+1}$ contains $F_j$, we may suppose that $g_{r+1}$ is the
blowing-up along $y_1=\cdots=y_a=y_{m+1}=0$, where $0 < a \leq m$.
The relevant chart is determined by the transformation
\[
(z_1,\ldots,z_n) \mapsto
(z_1z_{m+1},\ldots,z_az_{m+1},z_{a+1},\ldots,z_m,z_{m+1},\ldots,z_n).
\]
Because
\[
f \circ g_1 \circ \cdots \circ g_{r+1} = g^*_{r+1} \alpha \left(
\prod_{i=1}^{m} z_i^{N_i} \right) z_{m+1}^{\sum_{i=1}^{a} N_i}
\]
and
\[
(g_1 \circ \cdots \circ g_{r+1})^* dx = g^*_{r+1} \beta
\left(\prod_{i=1}^{m} z_i^{\nu_i-1} \right)
z_{m+1}^{\sum_{i=1}^{a} \nu_i} dz,
\]
we have to prove that
$N_{j,r}s_0+\nu_{j,r}=(N_{j,r}+\sum_{i=1}^{a}N_i)s_0+ (\nu_{j,r}+
\sum_{i=1}^{a} \nu_i)$. This follows from the fact that
$N_is_0+\nu_i=0$ for $i \in \{1,\ldots,a\}$. Remark that it
follows also from these calculations that
\[
N_{j,r} \equiv N_j \mbox{ mod } \mathrm{gcd}(N_1,\ldots,N_m).
\]

\vspace{0,5cm}

\noindent \textbf{(2.8)} \textsl{Example.} We give an illustration
which is easy and well known. Let $f=x_1^2+x_2^2$. Let
$X=\mathbb{Z}_p \times \mathbb{Z}_p$. We want to determine the
poles of Igusa's $p$-adic zeta function associated to $f$. Notice
that $-1$ is a square in $\mathbb{Q}_p$ if and only if $-1$ is a
square $\mathbb{Z}/(p)$ and $p \not= 2$.

If $-1$ is a square in $\mathbb{Q}_p$, then $(f,dx)$ has already
normal crossings. We obtain a good chart for $(f,dx)$ by applying
the coordinate transformation $(y_1,y_2) \mapsto
((y_1+y_2)/2,(y_1-y_2)/(2a))$, where $a$ denotes a square root of
$-1$. Because $|a|=1$ we obtain
\[
\int_{\mathbb{Z}_p \times \mathbb{Z}_p} |x_1^2+x_2^2|^s |dx_1
\wedge dx_2| = \int_{\mathbb{Z}_p \times \mathbb{Z}_p}
|y_1y_2|^{s} |dy_1 \wedge dy_2|.
\]
Consequently, the only candidate poles of $Z_f(s)$ are $-1+(2k\pi
\sqrt{-1})/(\log p)$, $k \in \mathbb{Z}$. They are all poles
because $b_{-2}=((p-1)/(p \log p))^2$ for each candidate pole.

If $-1$ is not a square in $\mathbb{Q}_p$, then $(f,dx)$ does not
have normal crossings at the origin. We obtain an embedded
resolution after one blowing-up $g$. Remark that the zero locus of
$f$ contains only the origin and that the zero locus of $f \circ
g$ is equal to the exceptional curve $E$ of $g$. We will use the
two charts on the blowing-up determined by $(y_1,y_2) \mapsto
(y_1y_2,y_2)$ and $(z_1,z_2) \mapsto (z_1,z_1z_2)$. The sets
$\{(y_1,y_2) \mid y_1 \in \mathbb{Z}_p \, , \, y_2=0 \}$ and
$\{(z_1,z_2) \mid z_1=0 \, , \, z_2 \in p \mathbb{Z}_p \}$ form a
partition of $E$. The candidate poles of $Z_f(s)$ are
$s_k=-1+(2k\pi \sqrt{-1})/(2 \log p)$, $k \in \mathbb{Z}$, and
each $b_{-1}$ is equal to
\[
\left( \frac{p-1}{2p \log p} \right) \left( \left[
\int_{\mathbb{Z}_p} |y_1^2+1|^s |dy_1|
\right]^{\mathrm{mc}}_{s=s_k} + \left[ \int_{p\mathbb{Z}_p}
|1+z_2^2|^s |dz_2| \right]^{\mathrm{mc}}_{s=s_k} \right).
\]
If $p \not= 2$, we have that $|1+x^2|=1$ for every $x \in
\mathbb{Z}_p$, so that $b_{-1}= (p^2-1)/(2p^2 \log p)$. If $p=2$,
we have that $|1+x^2|=1$ for every $x \in 2\mathbb{Z}_2$ and
$|1+x^2|=1/2$ for every $x \in 1+2\mathbb{Z}_2$, so that
$b_{-1}=1/(2 \log 2)$ if $k$ is even and $b_{-1}=0$ if $k$ is odd.

Remark that Igusa's $p$-adic zeta function of $x_1^2+x_2^2$ can be
calculated completely elementarily in all the cases.

\section{The vanishing results}
\subsection{Curves}
Let $X$ be an open and compact subset of $K^2$. Let $f$ be a
$K$-analytic function on $X$. Let $g:Y \rightarrow X$ be an
embedded resolution of $f$. Write $g=g_1 \circ \cdots \circ g_t:
Y=Y_t \rightarrow X=Y_0$ as a composition of blowing-ups $g_i: Y_i
\rightarrow Y_{i-1}$, $i \in \{1,\ldots,t\}$. The exceptional
curve of $g_i$ and also the strict transforms of this curve are
denoted by $E_i$. Let $\chi$ be a character of $R^{\times}$.

\vspace{0,5cm}

\noindent \textbf{Proposition.} \textsl{Let $r \in \{1,\ldots,t\}$
and let $P \in Y_{r-1}$ be the centre of the blowing-up $g_r$.
Suppose that the expected order of a candidate pole $s_0$
associated to $E_r$ is one. Suppose that there exists a chart
$(V,y=(y_1,y_2))$ centred at $P$ on which $f \circ g_1 \circ
\cdots \circ g_{r-1}$ is given by a power series with lowest
degree part a (non-constant) monomial. Then the contribution of
$E_r$ to the residue $b_{-1}$ of $Z_{f,\chi}(s)$ at $s_0$ is
zero.}

\vspace{0,2cm}

\noindent \textsl{Remark.} This proposition is essentially well
known. Our proof differs slightly from the ones in
\cite{Igusacomplexpowers} and \cite{Loeser1} because we will
calculate the contribution of $E_r$ to $b_{-1}$ just after the
creation of $E_r$ instead of on the embedded resolution. We
incorporate this proof here because the same technique will be
used in the proof of the more difficult result of section 3.2.

\vspace{0,2cm}

\noindent \textsl{Proof.} We may suppose that $(V,y)$ is a chart
centred at $P$ such that $f \circ g_1 \circ \cdots \circ g_{r-1} =
ey_1^k y_2^l + \theta$ and $(g_1 \circ \cdots \circ g_{r-1})^*dx =
\rho y_1^{c-1} y_2^{d-1} dy$ with $k,l \in \mathbb{Z}_{\geq 0}$,
$c,d \in \mathbb{Z}_{>0}$, $e \in K^{\times}$ and $\rho,\theta$
$K$-analytic functions satisfying $\rho(0,0) \not= 0$ and
$\mbox{mult}(\theta)>k+l$. We consider here the case that $k$ and
$l$ are both not zero. The case that $k$ or $l$ is zero can be
treated analogously.

We look at the chart $(O,z=(z_1,z_2))$ on $Y_r$ for which
$g_r(z_1,z_2)=(z_1,z_1z_2)$. Then
\begin{eqnarray*}
f \circ g_1 \circ \cdots \circ g_r & = & z_1^{k+l}
\left(ez_2^l+z_1 \frac{\theta(z_1,z_1z_2)}{z_1^{k+l+1}} \right)
\qquad \mbox{and}
\\ (g_1 \circ \cdots \circ g_r)^*dx & = & \rho(z_1,z_1z_2) z_1^{c+d-1} z_2^{d-1} dz.
\end{eqnarray*}
Remark that the equation of $E_r$ is $z_1=0$, that $N_r=k+l$ and
that $\nu_r=c+d$. Using the notation of (2.7), let $T_r=\{1,2\}$
and let $F_1$ be the origin of this chart. The contribution to
$b_{-1}$ of an open and compact subset $A$ of $E_r$ which is
contained in $O$ is equal to
\[
\left( \frac{q-1}{qN_r \log q} \right) \left[ \int_A
\chi(\mathrm{ac} \, e) \chi^l(\mathrm{ac} \, z_2) |e|^s
|\rho(0,0)| |z_2|^{ls+d-1} |dz_2| \right]^{\mathrm{mc}}_{s=s_0}.
\]

Let $(O',z'=(z_1',z_2'))$ be the chart on $Y_r$ for which
$g_r(z_1',z_2')=(z_1'z_2',z_2')$. The origin of this chart is the
point $F_2$. Analogously as before, we obtain that the
contribution to $b_{-1}$ of an open and compact subset $B$ of
$E_r$ which is contained in $O'$ is equal to
\[
\left( \frac{q-1}{qN_r \log q} \right) \left[ \int_B
\chi(\mathrm{ac} \, e) \chi^k(\mathrm{ac} \, z_1') |e|^s
|\rho(0,0)| |z_1'|^{ks+c-1} |dz_1'| \right]^{\mathrm{mc}}_{s=s_0}.
\]

Because $\chi^{N_r}=1$ (otherwise there are no candidate poles
associated to $E_r$) and because $k+l=N_r$, we have that
$\chi^k=1$ if and only if $\chi^l=1$.

\textsl{Case 1:} $\chi^k = \chi^l=1$. Then the contribution of
$E_r$ to $b_{-1}$ is equal to
\begin{eqnarray*}
\lefteqn{ \left( \frac{\chi(\mathrm{ac} \, e) |e|^{s_0}
|\rho(0,0)|(q-1)}{qN_r \log q} \right) \left( \left[ \int_R
|z_2|^{ls+d-1} |dz_2| \right]^{\mathrm{mc}}_{s=s_0} + \left[
\int_P |z_1'|^{ks+c-1} |dz_1'| \right]^{\mathrm{mc}}_{s=s_0}
\right) } \\ & = & \left( \frac{\chi(\mathrm{ac} \, e) |e|^{s_0}
|\rho(0,0)|(q-1)}{qN_r \log q} \right) \left( \frac{q-1}{q}
\frac{1}{1-q^{-\alpha_1}} + \frac{q-1}{q}
\frac{q^{-\alpha_2}}{1-q^{-\alpha_2}} \right) \\ & = & \left(
\frac{\chi(\mathrm{ac} \, e) |e|^{s_0} |\rho(0,0)|(q-1)}{qN_r \log
q} \right) \left( \frac{q-1}{q} \right) \left(
\frac{1-q^{-\alpha_2} +q^{-\alpha_2}-q^{-\alpha_1-\alpha_2}}
{(1-q^{-\alpha_1})(1-q^{-\alpha_2})} \right) \\ & = & 0.
\end{eqnarray*}
The last equality follows from $\alpha_1 + \alpha_2=0$, which is
relation 2 of (2.7).

\textsl{Case 2:} $\chi^k \not= 1$ and $\chi^l \not= 1$. Then the
contribution of $E_r$ to $b_{-1}$ is equal to zero because both
terms in the sum
\[
\left[ \int_R \chi^l(\mathrm{ac} \, z_2) |z_2|^{ls+d-1} |dz_2|
\right]^{\mathrm{mc}}_{s=s_0} + \left[ \int_P \chi^k(\mathrm{ac}
\, z_1') |z_1'|^{ks+c-1} |dz_1'| \right]^{\mathrm{mc}}_{s=s_0}
\]
are equal to zero. $\qquad \Box$

\subsection{Surfaces}
Let $X$ be an open and compact subset of $K^3$. Let $f$ be a
$K$-analytic function on $X$. Let $g:Y = Y_t \rightarrow X = Y_0$
be an embedded resolution of $f$ which is a composition $g_1 \circ
\cdots \circ g_t$ of blowing-ups $g_i:Y_i \rightarrow Y_{i-1}$
with centre a $K$-analytic closed submanifold which has only
normal crossings with the union of the exceptional surfaces in
$Y_{i-1}$ and with exceptional surface $E_i$.

\vspace{0,5cm}

\noindent \textbf{Proposition.} \textsl{Let $r \in \{1,\ldots,t\}$
and let $P \in Y_{r-1}$ be the centre of the blowing-up $g_r$.
Suppose that the expected order of a candidate pole $s_0$
associated to $E_r$ is one. Suppose that there exists a chart
$(V,y=(y_1,y_2,y_3))$ centred at $P$ on which $f \circ g_1 \circ
\cdots \circ g_{r-1}$ is given by a power series with lowest
degree part of the form $ey_1^ky_2^ly_3^m(y_1+y_2)^n$, with $e \in
K^{\times}$ and $k,l,m,n \in \mathbb{Z}_{\geq 0}$. Then the
contribution of $E_r$ to the residue $b_{-1}$ of $Z_f(s)$ at $s_0$
is zero.}

\vspace{0,2cm}

\noindent \textsl{Proof.} We may suppose that $f \circ g_1 \circ
\cdots \circ g_{r-1} = ey_1^k y_2^ly_3^m(y_1+y_2)^n + \theta$ and
$(g_1 \circ \cdots \circ g_{r-1})^*dx = \rho y_1^{a-1} y_2^{b-1}
y_3^{c-1} (y_1+y_2)^{d-1} dy$ with $a,b,c,d \in \mathbb{Z}_{>0}$
and $\rho,\theta$ $K$-analytic functions satisfying $\rho(0,0)
\not= 0$ and $\mbox{mult}(\theta)>k+l+m+n$. Remark that at least
one of the numbers $a,b,d$ is equal to $1$. We consider here the
case that $k,l,m$ and $n$ are all different from zero. The other
cases are treated analogously. Let $T_r=\{ 1,2,3,4 \}$ and suppose
that $F_i$, $i \in \{ 1,2,3 \}$, is given by $y_i=0$ and that
$F_4$ is given by $y_1+y_2=0$ in the homogeneous coordinates
$(y_1:y_2:y_3)$ on $E_r \subset Y_r$.

Analogously as in Section 3.1, we can calculate $f \circ g_1 \circ
\cdots \circ g_r$ and $(g_1 \circ \cdots \circ g_r)^*dx$ in the
three charts on $Y_r$ for which respectively
$g_r(z_1,z_2,z_3)=(z_1z_3,z_2z_3,z_3)$,
$g_r(z_1',z_2',z_3')=(z_1'z_2',z_2',z_2'z_3')$ and
$g_r(z_1'',z_2'',z_3'')=(z_1'',z_1''z_2'',z_1''z_3'')$. The
contribution of $E_r$ to the residue $b_{-1}$ of $Z_f(s)$ at $s_0$
turns out to be $\kappa |e|^{s_0} |\rho(0,0,0)|$ times
\vspace{-0,1cm}
\begin{picture}(150,0)(0,0)
\put(135,-23){($*$)}
\end{picture}
\begin{eqnarray*}
\lefteqn{ \left[ \int_{P \times P} |z_1|^{ks+a-1} |z_2|^{ls+b-1}
|z_1+z_2|^{ns+d-1} |dz_1 \wedge dz_2|
\right]^{\mathrm{mc}}_{s=s_0} } \\
 & + & \left[ \int_{P \times R} |z_1'|^{ks+a-1} |z_3'|^{ms+c-1} |z_1'+1|^{ns+d-1}
|dz_1' \wedge dz_3'| \right]^{\mathrm{mc}}_{s=s_0} \\ & + & \left[
\int_{R \times R} |z_2''|^{ls+b-1} |z_3''|^{ms+c-1}
|1+z_2''|^{ns+d-1} |dz_2'' \wedge dz_3''|
\right]^{\mathrm{mc}}_{s=s_0}.
\end{eqnarray*}
Consequently, we have to prove that this expression is equal to
zero.

To calculate the first term in $(*)$, we partition $P \times P$
into
\begin{eqnarray*}
A_1 & = & \{ (z_1,z_2) \in P \times P \mid \mathrm{ord} \, z_1 >
\mathrm{ord} \, z_2 \} = \bigsqcup_{i \in \mathbb{Z}_{>0}} \{
(z_1,z_2) \mid \mathrm{ord} \, z_1 > \mathrm{ord} \, z_2=i \} \\
A_2 & = & \{ (z_1,z_2) \in P \times P \mid \mathrm{ord} \, z_1 <
\mathrm{ord} \, z_2 \} = \bigsqcup_{i \in \mathbb{Z}_{>0}} \{
(z_1,z_2) \mid i= \mathrm{ord} \, z_1 < \mathrm{ord} \, z_2 \} \\
A_3 & = & \{ (z_1,z_2) \in P \times P \mid \mathrm{ord} \, z_1 =
\mathrm{ord} \, z_2 \} = \bigsqcup_{i \in \mathbb{Z}_{>0}} \{
(z_1,z_2) \mid \mathrm{ord} \, z_1 = \mathrm{ord} \, z_2=i \}
\end{eqnarray*}
The contribution of $A_1$ to the first term in $(*)$ is equal to
\begin{eqnarray}
\lefteqn{ \left[ \sum_{i=1}^{\infty} \int_{P^{i+1}} \left(
\int_{P^i \setminus P^{i+1}} |z_1|^{ks+a-1} |z_2|^{ls+b-1}
|z_1+z_2|^{ns+d-1} |dz_2| \right) |dz_1|
\right]^{\mathrm{mc}}_{s=s_0} } \nonumber \\ & = & \left[
\sum_{i=1}^{\infty} \frac{q-1}{q} q^{-i} q^{-i(ls+b-1)}
q^{-i(ns+d-1)} \int_{P^{i+1}} |z_1|^{ks+a-1} |dz_1|
\right]^{\mathrm{mc}}_{s=s_0} \nonumber \\ & = & \left[ \left(
\frac{q-1}{q} \right)^2 \frac{1}{q^{ks+a}-1} \sum_{i=1}^{\infty}
q^{-i(ks+a+ls+b+ns+d-1)} \right]^{\mathrm{mc}}_{s=s_0} \nonumber
\\ & = & \left[ \left( \frac{q-1}{q} \right)^2
\frac{1}{(q^{ks+a}-1)(q^{ks+a+ls+b+ns+d-1}-1)} \right]^{\mathrm{mc}}_{s=s_0}
\nonumber \\ & = & \left( \frac{q-1}{q} \right)^2
\frac{1}{(q^{\alpha_1}-1) (q^{\alpha_1+\alpha_2+\alpha_4-1}-1)}.
\end{eqnarray}
Analogously, we obtain that the contribution of $A_2$ to the first
term in $(*)$ is equal to
\begin{eqnarray}
\left( \frac{q-1}{q} \right)^2
\frac{1}{(q^{\alpha_2}-1)(q^{\alpha_1+\alpha_2+\alpha_4-1}-1)}.
\end{eqnarray}
The contribution of $A_3$ to the first term in $(*)$ is
\[
\left[ \sum_{i=1}^{\infty} \int_{(P^i \setminus P^{i+1})^2}
|z_1|^{ks+a-1} |z_2|^{ls+b-1} |z_1+z_2|^{ns+d-1} |dz_1 \wedge
dz_2| \right]^{\mathrm{mc}}_{s=s_0}.
\]
One can verify (see \cite[Section 3.3.2]{Segers}) that this is
equal to
\begin{eqnarray}
\left( \frac{q-1}{q} \right)^2
\frac{1}{(q^{\alpha_4}-1)(q^{\alpha_1+\alpha_2+\alpha_4-1}-1)}
\hspace{1cm} \\  + \left( \frac{q-1}{q} \right) \left(
\frac{q-2}{q} \right)
\frac{1}{q^{\alpha_1+\alpha_2+\alpha_4-1}-1}.
\end{eqnarray}
The second term of $(*)$ is equal to
\begin{eqnarray}
\left[ \int_{P} |z_1|^{ks+a-1} |dz_1| \int_R |z_3|^{ms+c-1} |dz_3|
\right]^{\mathrm{mc}}_{s=s_0} = \left( \frac{q-1}{q} \right)^2
\frac{1}{(q^{\alpha_1}-1)(1-q^{-\alpha_3})}.
\end{eqnarray}
The third term of $(*)$ is equal to
\begin{eqnarray}
\lefteqn{ \left[ \int_R |z_2|^{ls+b-1} |1+z_2|^{ns+d-1} |dz_2|
\int_R |z_3|^{ms+c-1} |dz_3| \right]^{\mathrm{mc}}_{s=s_0} }
\nonumber
\\ & = & \left[ \left( \int_{R \setminus (P \cup -1+P)}
|z_2|^{ls+b-1} |1+z_2|^{ns+d-1} |dz_2| + \int_P |z_2|^{ls+b-1}
|1+z_2|^{ns+d-1} |dz_2| \right. \right. \nonumber \\ & & \quad +
\left. \left. \int_{-1+P} |z_2|^{ls+b-1} |1+z_2|^{ns+d-1} |dz_2|
\right) \left( \int_R |z_3|^{ms+c-1} |dz_3| \right)
\right]^{\mathrm{mc}}_{s=s_0} \nonumber \\ & = & \left[ \left( 1-
\frac{2}{q} + \int_P |z_2|^{ls+b-1} |dz_2| + \int_{-1+P}
|1+z_2|^{ns+d-1} |dz_2| \right) \left( \int_R |z_3|^{ms+c-1}
|dz_3| \right) \right]^{\mathrm{mc}}_{s=s_0} \nonumber \\ & = &
\left(1- \frac{2}{q} + \frac{q-1}{q} \frac{1}{q^{\alpha_2}-1} +
\frac{q-1}{q} \frac{1}{q^{\alpha_4}-1} \right) \left(
\frac{q-1}{q} \frac{1}{1-q^{-\alpha_3}} \right) \nonumber \\ & = &
\left( \frac{q-1}{q} \right) \left( \frac{q-2}{q} \right)
\frac{1}{1-q^{-\alpha_3}}
\\ & & \quad + \left( \frac{q-1}{q} \right)^2
\frac{1}{(q^{\alpha_2}-1)(1-q^{-\alpha_3})} \\ &  & \quad + \left(
\frac{q-1}{q} \right)^2
\frac{1}{(q^{\alpha_4}-1)(1-q^{-\alpha_3})}.
\end{eqnarray}
Relation 2 of (2.7) is $\alpha_1+\alpha_2+\alpha_3+\alpha_4-1=0$,
so that we obtain
\begin{eqnarray*}
\frac{1}{q^{\alpha_1+\alpha_2+\alpha_4-1}-1} +
\frac{1}{1-q^{-\alpha_3}} & = &
\frac{1-q^{-\alpha_3}+q^{\alpha_1+\alpha_2+\alpha_4-1}-1}
{(q^{\alpha_1+\alpha_2+\alpha_4-1}-1)(1-q^{-\alpha_3})} \\ & = &
\frac{q^{\alpha_1+\alpha_2+\alpha_3+\alpha_4-1}-1}
{(q^{\alpha_1+\alpha_2+\alpha_4-1}-1)(q^{\alpha_3}-1)} \\ & = & 0,
\end{eqnarray*}
and consequently $(3)+(7)=0$. Analogously, we obtain that
$(4)+(9)=(5)+(10)=(6)+(8)=0$. Consequently, the contribution of
$E_r$ to $b_{-1}$ is equal to zero. $\qquad \Box$

\vspace{0,2cm}

\noindent \textsl{Remark.} Let $\chi$ be a character of
$R^{\times}$. Suppose that we are in the analogous situation to
this proposition with $Z_{f,\chi}(s)$. Then one can show that the
contribution of $E_r$ to the residue $b_{-1}$ of $Z_{f,\chi}(s)$
at $s_0$ is zero. The proof consists of very long calculations
involving character sums and is written down in \cite[Chapter
4]{Segers}.

\section{Determination of the smallest poles}

The main ideas and results of this section have the same flavour
as those in \cite{SegersVeys}, where the local topological zeta
function is studied. However here the situation is more
complicated because the field $K$ is not algebraically closed.

\subsection{Curves}
\noindent \textbf{(4.1.1)} In this section we will determine
$\mathcal{P}_2^K \cap ]-\infty,-1/2[$. Let $f$ be a $K$-analytic
function on an open and compact subset of $K^2$ and let $g$ be the
\textsl{minimal} embedded resolution of $f$. The poles of $Z_f(s)$
with real part less than $-1/2$ and different from $-1$ are only
associated to exceptional curves. Consequently, these poles are
completely determined by the germs of $f$ at the points where $f$
does not have normal crossings. It is thus sufficient to study the
germs of $K$-analytic functions at the origin, which will be
identified with the convergent power series. The set of all
convergent power series in the variables $x$ and $y$ is
classically denoted by $K\!\!<\!\!<\!\!x,y\!\!>\!\!>$.

\vspace{0,5cm}

\noindent \textbf{(4.1.2)} Let $f \in
K\!\!<\!\!<\!\!x,y\!\!>\!\!>$. Let $g:Y \rightarrow X$ be the
\textsl{minimal} embedded resolution of a representative of $f$.
Write $g=g_1 \circ \cdots \circ g_t: Y=Y_t \rightarrow X=Y_0$ as a
composition of blowing-ups $g_i: Y_i \rightarrow Y_{i-1}$, $i \in
\{1,\ldots,t\}$. The exceptional curve of $g_i$ and also the
strict transforms of this curve are denoted by $E_i$. Let $T$ be
as in (2.1) and obviously we suppose that $\{1,\ldots,t\} \subset
T$.

Let $k \in \{1,\ldots,t\}$. Let $P \in Y_k$ be a point on an
exceptional curve, i.e., a point which is mapped to the origin
under the map $g_1 \circ \cdots \circ g_k$. The strict transform
of $f$ around $P$ is defined as the germ at $P$ of the
$K$-analytic function $f \circ g_1 \circ \cdots \circ g_k$ divided
by the highest possible powers of local equations of exceptional
curves through $P$. Remark that the strict transform of $f$ around
$P$ is defined modulo the germ of a $K$-analytic function which
does not vanish at $P$ as a factor.

We call a complex number `a candidate pole of $Z_f(s)$' if it is a
candidate pole associated to an $E_i$, $i \in T$, satisfying $0
\in g(E_i)$. A candidate pole of $Z_f(s)$ is called a pole of
$Z_f(s)$ if there exists an arbitrarily small neighbourhood of $0$
for which it is a pole.

The following lemma is trivial.

\vspace{0,5cm}

\noindent \textbf{(4.1.3) Lemma.} \textsl{Suppose that we have
blown up $k$ times but we do not yet have an embedded resolution.
Let $P$ be a point at which $f \circ g_1 \circ \cdots \circ g_k$
does not have normal crossings. Let $\mu$ be the multiplicity at
$P$ of the strict transform of $f$ around $P$ and let $g_{k+1}$ be
the blowing-up at $P$.}

\noindent \hspace{0,5cm} \textsl{(a) Suppose that two exceptional
curves $E_i$ and $E_j$ contain $P$. Then $-\nu_{k+1}/N_{k+1}$ is
equal to $-(\nu_i+\nu_j)/(N_i+N_j+\mu)$ and this value is larger
than $\min\{-\nu_i/N_i,-\nu_j/N_j\}$.}

\noindent \hspace{0,5cm} \textsl{(b) Suppose that exactly one
exceptional curve $E_i$ contains $P$ and that $\mu \geq 2$. Then
$E_{k+1}$ has numerical data $(N_i+\mu,\nu_i+1)$ and
$-(\nu_i+1)/(N_i+\mu)$ lies between $-1/\mu$ and $-\nu_i/N_i$.}

\noindent \hspace{0,5cm} \textsl{(c) Suppose that exactly one
exceptional curve $E_i$ contains $P$ and that $\mu=1$. Remark that
the two curves are tangent at $P$ because we do not have normal
crossings at $P$. Let $g_{k+2}$ be the blowing-up at $E_i \cap
E_{k+1}$. Remark that we do not have to blow up at a point of
$E_{k+1}$ anymore. The numerical data of $E_{k+2}$ are
$(2N_i+2,2\nu_i+1)$, and $-(2\nu_i+1)/(2N_i+2)$ lies between
$-1/2$ and $-\nu_i/N_i$. Let $s_0$ be a candidate pole associated
to $E_{k+1}$. Because $s_0$ is not a candidate pole associated to
$E_{k+2}$, which is a consequence of $-\nu_{k+1}/N_{k+1} \not=
-\nu_{k+2}/N_{k+2}$, the contribution of $E_{k+1}$ to the
coefficient $b_{-2}$ in the Laurent series of $Z_f(s)$ at $s_0$ is
zero. It follows from the proposition in $3.1$ that $E_{k+1}$ does
not give a contribution to the residue $b_{-1}$ of $Z_f(s)$ at
$s_0$.}

\vspace{0,5cm}

\noindent \textbf{(4.1.4)} Suppose that after some blowing-ups,
the pullback of $f$ does not have normal crossings at a point $P$.
Suppose also that the real parts of the candidate poles associated
to the exceptional curves through $P$ are all larger than or equal
to $-1/2$. Then it follows from the above lemma that the
components above $P$ in the final resolution do not give a
contribution to a candidate pole with real part less than $-1/2$.

\vspace{0,2cm} \noindent \textbf{Corollary.} \textsl{Zeta
functions of convergent power series of multiplicity at least $4$
do not have a pole with real part in} $]-\infty,-1/2[ \setminus
\{-1\}$.

\vspace{0,2cm} \noindent Indeed, every exceptional curve in the
minimal embedded resolution of $f$ lies above a point of $E_1$
(considered in the stage when it is created), which has a
candidate pole with real part larger than or equal to $-1/2$.

\vspace{0,5cm}

\noindent \textbf{(4.1.5)} To deal with multiplicity $2$ and $3$,
we will study an `easier' element of
$K\!\!<\!\!<\!\!x,y\!\!>\!\!>$. We will use the following theorem
(see \cite[Theorem 2.3.1]{Igusabook}).

\vspace{0,2cm} \noindent \textsc{Weierstrass Preparation Theorem.}
\\ If $f(z_1,\ldots,z_{n-1},w)=f(z,w) \in
K\!\!<\!\!<\!\!z,w\!\!>\!\!>$ is not identically zero on the
$w$-axis, then $f$ can be written uniquely as
$f=(w^e+a_1(z)w^{e-1}+\cdots+a_e(z))h$, where $a_i(z) \in
K\!\!<\!\!<\!\!z\!\!>\!\!>$ satisfies $a_i(0)=0$ and $h \in
K\!\!<\!\!<\!\!z,w\!\!>\!\!>$ satisfies $h(0) \not= 0$.

\vspace{0,2cm} \noindent Because $h(0) \not= 0$ implies that $|h|$
is constant on a neighbourhood of  $0$, we have that Igusa's
$p$-adic zeta functions of $f$ and
$w^e+a_1(z)w^{e-1}+\cdots+a_e(z)$ have the same poles. After an
appropriate coordinate transformation, the desired form will
appear. For example, the coordinate transformation $(z,w) \mapsto
(z,w-a_1(z)/e)$ cancels the term $a_1(z)w^{e-1}$.

\vspace{0,5cm}

\noindent \textbf{(4.1.6)} \textsl{Example.} Let $f \in
K\!\!<\!\!<\!\!x,y\!\!>\!\!>$ have multiplicity $3$ and let
$f_3=y^3+xy^2=y^2(y+x)$ be the homogeneous part of $f$ of degree
$3$. By the Weierstrass preparation theorem, we may work with a
function of the form $y^3+a_1(x)y^2+a_2(x)y+a_3(x)$, with
$\mbox{mult}(a_1(x))=1$, $\mbox{mult}(a_2(x)) \geq 3$ and
$\mbox{mult}(a_3(x)) \geq 4$. One can check that there exists a
coordinate transformation $(x,y) \mapsto (x,y-k(x))$ such that the
function becomes of the form $y^3+b_1(x)y^2+b_3(x)$, with
$\mbox{mult}(b_1(x))=1$ and $\mbox{mult}(b_3(x)) \geq 4$. After
another coordinate transformation, we get the form
$y^3+xy^2+g(x)$, with $\mbox{mult}(g(x)) \geq 4$.

\vspace{0,5cm}

\noindent \textbf{(4.1.7) Theorem.} \textsl{We have}
\[ \mathcal{P}_2^K \cap \left]-\infty,-\frac{1}{2}\right[ = \left\{ \left.
-\frac{1}{2}-\frac{1}{i} \right| i \in \mathbb{Z}_{>1} \right\} \]
\textsl{and at most one value in $]-1,-1/2]$ is the real part of a
pole of a fixed Igusa's $p$-adic zeta function. Moreover, if $f
\in K\!\!<\!\!<\!\!x,y\!\!>\!\!>$ has multiplicity at least $4$,
then $Z_f(s)$ has no pole with real part in $]-\infty,-1/2[
\setminus \{-1\}$.}

\vspace{0,2cm} \noindent \emph{Proof.} Because the calculations
are analogous to the calculations in \cite{SegersVeys} for the
local topological zeta function, we do not treat all the cases in
this paper.

\vspace{0,4cm} \noindent \hspace{1cm} (a) Suppose that $f$ is an
element of $K\!\!<\!\!<\!\!x_1,x_2\!\!>\!\!>$ with multiplicity 2.
When we apply the ideas of (4.1.5), we see that it is enough to
consider $x_1^2$ and $x_1^2+ax_2^l$, with $l \in \mathbb{Z}_{>1}$
and $a \in K^{\times}$. If $f=x_1^2$, the candidate poles of
$Z_f(s)$ are $-1/2+(k \pi \sqrt{-1})/(\log q)$, $k \in
\mathbb{Z}$. If $l=2$, the calculations are analogous as in (2.8).
If $l$ is odd, write $l=2r+1$. After $r$ blowing-ups, the strict
transform of $f^{-1}\{0\}$ is non-singular and tangent to $E_r$.
The numerical data of $E_i$, $i=1,\ldots,r$, are $(2i,i+1)$. To
get the minimal embedded resolution, we now blow up twice. Let
$E_0$ be the strict transform of $f^{-1}\{0\}$. Remark that
$T=\{0,1,\ldots,r+2\}$. The dual resolution graph and the
numerical data are given below.
\\
\begin{picture}(150,20)(-5,2)
\put(0,15){\line(1,0){25}} \put(28,15){$\ldots$}
\put(35,15){\line(1,0){25}} \put(0,15){\circle*{1.5}}
\put(10,15){\circle*{1.5}} \put(20,15){\circle*{1.5}}
\put(40,15){\circle*{1.5}} \put(50,15){\circle*{1.5}}
\put(60,15){\circle*{1.5}} \put(50,14.5){\line(0,-1){8.75}}
\put(50,5){\circle{1.5}} \put(-2,17){$E_1$} \put(8,17){$E_2$}
\put(18,17){$E_3$} \put(38,17){$E_r$} \put(48,17){$E_{r+2}$}
\put(58,17){$E_{r+1}$} \put(52,5){$E_0$}
\put(75,5){\shortstack[l]{$E_1(2,2)$
\\ $E_2(4,3)$ \\ $E_3(6,4)$}}
\put(100,5){\shortstack[l]{$E_r(2r,r+1)$
\\ $E_{r+1}(2r+1,r+2)$ \\ $E_{r+2}(4r+2,2r+3)$}}
\end{picture}
It follows from section 3.1 that the candidate poles associated to
$E_1,\ldots,E_{r+1}$ are not poles. The other candidate poles have
real part $-1$ or $-(2r+3)/(4r+2)=-1/2-1/(2r+1)$. We calculate the
residue of $Z_f(s)$ at the candidate pole $s_0=-1/2-1/(2r+1)$.
Because
\begin{eqnarray*}
\lefteqn{ \left[ \int_{aR} |y_1|^{(2r+1)s+r+1} |y_1+a|^s |dy_1|
\right]^{\mathrm{mc}}_{s=s_0} } \\ & = & |a|^{-1/(2r+1)} \left[
\int_{R} |y|^{(2r+1)s+r+1} |y+1|^s |dy|
\right]^{\mathrm{mc}}_{s=s_0}  \\ & = & |a|^{-1/(2r+1)} \left[
\int_{R \setminus (-1+P)} |y|^{(2r+1)s+r+1} |dy| + \int_{-1+P}
|y+1|^s |dy| \right]^{\mathrm{mc}}_{s=s_0}
\\ & = & |a|^{-1/(2r+1)} \left( \frac{q-2}{q} + \frac{q-1}{q}
\frac{1}{q^{\alpha_{r+1}}-1} + \frac{q-1}{q}
\frac{1}{q^{\alpha_0}-1} \right)
\end{eqnarray*}
and
\begin{eqnarray*}
\left[ \int_{\frac{1}{a}P} |y_2|^{2rs+r} |1+ay_2|^s |dy_2|
\right]^{\mathrm{mc}}_{s=s_0} & = & |a|^{-1/(2r+1)} \left[
\int_{P} |y|^{2rs+r} |1+y|^s |dy| \right]^{\mathrm{mc}}_{s=s_0}
\\ & = & |a|^{-1/(2r+1)} \left[ \int_{P} |y|^{2rs+r} |dy|
\right]^{\mathrm{mc}}_{s=s_0}
\\ & = & |a|^{-1/(2r+1)} \frac{q-1}{q} \frac{1}{q^{\alpha_r}-1}
\end{eqnarray*}
the residue of $Z_f(s)$ at the candidate pole $s_0=-1/2-1/(2r+1)$
is
\[
|a|^{-1/(2r+1)} \left( \frac{q-2}{q} + \frac{q-1}{q}
\frac{1}{q^{\alpha_{r+1}}-1} + \frac{q-1}{q}
\frac{1}{q^{\alpha_0}-1} + \frac{q-1}{q} \frac{1}{q^{\alpha_r}-1}
\right)
\]
multiplied by the factor $\kappa$ which is different from zero
(see (2.6)). Because $\alpha_{r+1}=(2r+1)s_0+r+2=1/2
> 0$, $\alpha_0=s_0+1=1/2-1/(2r+1)
> 0$ and $\alpha_r=2rs_0+r+1=1/(2r+1)
> 0$, we have that the last three terms of this expression are
strictly positive. Consequently the whole expression is strictly
positive and thus different from zero, so that $-1/2-1/(2r+1)$ is
a pole of $Z_f(s)$.

If $l$ is even and larger than $2$, write $l=2r$. We have to blow
up $r$ times to obtain an embedded resolution. We have $E_1(2,2)$,
$E_2(4,3)$, $E_3(6,4)$, $\ldots$, $E_{r-1}(2r-2,r)$,
$E_r(2r,r+1)$. We obtain the first dual resolution graph if $-a$
is a square in $K$. Otherwise, we obtain the second dual
resolution graph.
\\
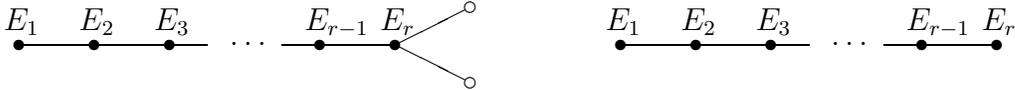
\begin{picture}(150,20)(-5,2)
\put(0,10){\line(1,0){25}} \put(28,10){$\ldots$}
\put(35,10){\line(1,0){15}} \put(0,10){\circle*{1.5}}
\put(10,10){\circle*{1.5}} \put(20,10){\circle*{1.5}}
\put(40,10){\circle*{1.5}} \put(50,10){\circle*{1.5}}
\put(50,10){\line(2,1){9.3}} \put(50,10){\line(2,-1){9.3}}
\put(60,15){\circle{1.5}} \put(60,5){\circle{1.5}}
\put(-2,12){$E_1$} \put(8,12){$E_2$} \put(18,12){$E_3$}
\put(38,12){$E_{r-1}$} \put(48,12){$E_r$}

\put(80,10){\line(1,0){25}} \put(108,10){$\ldots$}
\put(115,10){\line(1,0){15}} \put(80,10){\circle*{1.5}}
\put(90,10){\circle*{1.5}} \put(100,10){\circle*{1.5}}
\put(120,10){\circle*{1.5}} \put(130,10){\circle*{1.5}}
\put(78,12){$E_1$} \put(88,12){$E_2$} \put(98,12){$E_3$}
\put(118,12){$E_{r-1}$} \put(128,12){$E_r$}
\end{picture}
It follows from section 3.1 that the candidate poles associated to
$E_1,\ldots,E_{r-1}$ are not poles. The other candidate poles have
real part $-1$ or $-(r+1)/(2r)=-1/2-1/(2r)$ in the first case and
$-(r+1)/(2r)=-1/2-1/(2r)$ in the second case. Now we  prove that
$-1/2-1/(2r)$ is an element of $\mathcal{P}_2^K$. Suppose first
that $p \not= 2$. Then there exists an element $a$ of $K$ with
norm $1$ for which $-a$ is not a square in $K$. For such an $a$,
the residue of $Z_f(s)$ at $s_0=-1/2-1/(2r)$ is the non-zero
factor $\kappa$ times
\[
\frac{q-1}{q} \frac{1}{q^{\alpha_{r-1}}-1}+1.
\]
Suppose now that $p=2$. Remark that every element of the residue
field is a square in this case. Let $b \in R^{\times}$. If $b' \in
b+P$, then $b'^2-b^2 \in P^2$. Consequently, there exists an $a
\in -b^2+P$ such that $|a+x^2|=1/q$ for all $x \in b+P$. For such
an $a$, the residue of $Z_f(s)$ at $s_0=-1/2-1/(2r)$ is the
non-zero factor $\kappa$ times
\[
\frac{q-1}{q} \frac{1}{q^{\alpha_{r-1}}-1} + \frac{q-1}{q} +
\frac{1}{q} \left( \frac{1}{q} \right)^{s_0}.
\]
Because $\alpha_{r-1}=(2r-2)s_0+r=1/r > 0$, we obtain in the two
cases that this residue is strictly positive, which implies that
$-1/2-1/(2r)$ is a pole.

Our conclusion of part (a) is thus
\begin{eqnarray*} \{ s_0 \mid \exists f \in K\!\!<\!\!<\!\!x_1,x_2\!\!>\!\!> & : &
\mbox{mult}(f)=2 \mbox{ and } Z_f(s) \mbox{ has a pole with real
part } s_0 \}
\\ & = & \left\{ \left. -\frac{1}{2}-\frac{1}{i} \right| i \in
\mathbb{Z}_{>1} \right\} \cup \left\{ -\frac{1}{2} \right\}.
\end{eqnarray*}
Remark that Newton polyhedra could also be used to deal with (a),
see \cite{DenefHoornaert}.

\vspace{0,4cm} \noindent \hspace{1cm} (b) Suppose that $f$ is an
element of $K\!\!<\!\!<\!\!x_1,x_2\!\!>\!\!>$ with multiplicity 3.
Up to an affine coordinate transformation, there are three cases
for $f_3$.

\noindent \hspace{0,5cm} We consider the case that $f_3$ is a
product of three different linear factors over $K^{\mathrm{alg \,
cl}}$. Then we obtain an embedded resolution after one blowing-up.
There are three possibilities for the dual resolution graph,
depending on whether $f_3$ splits into linear factors over $K$,
$f_3$ is a product of a linear factor and an irreducible factor of
degree $2$ over $K$ or $f_3$ is irreducible over $K$. The dual
resolution graphs are respectively \\
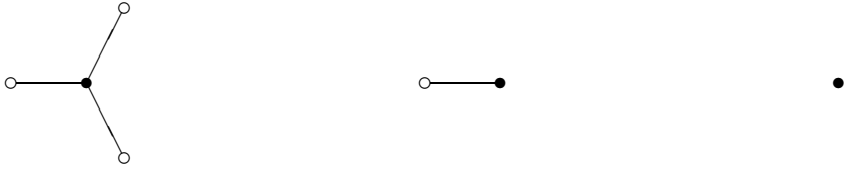
\begin{picture}(150,24)(-5,-2)
\put(20,10){\circle*{1.5}} \put(20,10){\line(-1,0){9.3}}
\put(20,10){\line(1,2){4.7}} \put(20,10){\line(1,-2){4.7}}
\put(10,10){\circle{1.4}} \put(25,20){\circle{1.4}}
\put(25,0){\circle{1.4}}

\put(75,10){\circle*{1.5}} \put(65,10){\circle{1.4}}
\put(75,10){\line(-1,0){9.3}}

\put(120,10){\circle*{1.5}}
\end{picture}
The equations of $f_3 \circ g$ in the charts determined by
$(y_1,y_2) \mapsto (y_1,y_1y_2)$ and $(z_1,z_2) \mapsto
(z_1z_2,z_2)$ are respectively of the form $y_1^3h_1$ and
$z_2^3h_2$. In the last case for example, we have that $h_1$ and
$h_2$ are non-vanishing on the exceptional curve. \\ The real
parts of the candidate poles of $Z_f(s)$ are $-1$ and
$-2/3=-1/2-1/6$ in the first two cases and $-2/3=-1/2-1/6$ in the
last case.

\noindent \hspace{0,5cm} The other cases are treated in
\cite{SegersVeys} for the topological zeta function and are very
similar for Igusa's $p$-adic zeta function.

\vspace{0,4cm} \noindent \hspace{1cm} (c) \hfill Suppose \hfill
that \hfill $f$ \hfill is \hfill an \hfill element \hfill of
\hfill $K\!\!<\!\!<\!\!x_1,x_2\!\!>\!\!>$ \hfill with \hfill
multiplicity \hfill at \hfill least \hfill 4. \hfill We \hfill
explained \hfill in \hfill (4.1.4) \hfill that \hfill $Z_f(s)$
\hfill has \hfill no \hfill pole \hfill with \hfill real \hfill
part \hfill in
\newline $]-\infty,-1/2[ \setminus \{-1\}$. $\qquad \Box$

\vspace{0,5cm}

\noindent \textbf{(4.1.8)} Let $\chi$ be a character of
$R^{\times}$. For $n \in \mathbb{Z}_{>0}$, we define the set
$\mathcal{P}_{n,\chi}^K$ by
\[
\mathcal{P}_{n,\chi}^K := \{ s_0 \mid \exists f \in F_n^K \, : \,
Z_{f,\chi}(s) \textsl{\mbox{ has a pole with real part }} s_0 \}.
\]

\vspace{0,2cm}

\noindent \textbf{Theorem.} \textsl{We have}
\[
\mathcal{P}_{2,\chi}^K \cap \left]-\infty,-\frac{1}{2}\right[
\subset \left\{ \left. -\frac{1}{2}-\frac{1}{i} \right| i \in
\mathbb{Z}_{>1} \right\} \] \textsl{and at most one value in
$]-1,-1/2]$ is the real part of a pole of a fixed Igusa's $p$-adic
zeta function.}

\vspace{0,2cm}

\noindent \emph{Proof.} This inculsion is proved in the same way
as the analogous inclusion of the previous theorem. Again we need
the proposition in 3.1. $\qquad \Box$

\subsection{Surfaces}
In this section, we prove the following theorem.

\vspace{0,2cm}

\noindent \textbf{(4.2.0) Theorem.} \textsl{We have}
\[ \mathcal{P}_3^K \cap ]-\infty,-1[ = \left\{ \left. -1-\frac{1}{i} \right| i \in
\mathbb{Z}_{>1} \right\}. \] \textsl{Moreover, if  $f \in
K\!\!<\!\!<\!\!x,y,z\!\!>\!\!>$ has multiplicity at least $3$,
then $Z_f(s)$ has no pole with real part less than $-1$.}

\vspace{0,2cm}

\noindent \textsl{Remark.} (i) It is a priori not obvious that the
smallest value of $\mathcal{P}_3^K$ is $-3/2$. This is in contrast
with the fact that it easily follows from lemma 4.1.3 that the
smallest value of $\mathcal{P}_2^K$ is $-1$. \\ (ii) Let $\chi$ be
a character of $R^{\times}$. Then one proves analogously as below
that an element of $\mathcal{P}_{3,\chi}^{K}$ less than $-1$ is of
the form $-1-1/i$, $i \in \mathbb{Z}_{>1}$. Using the remark in
section 3.2, the arguments below will also imply that
$Z_{f,\chi}(s)$ has no pole with real part less than $-1$ if $f
\in K\!\!<\!\!<\!\!x,y,z\!\!>\!\!>$ has multiplicity at least $3$.

\subsubsection{Multiplicity 2}
\noindent \textbf{(4.2.1.1)} Let $f(x)$, $x=(x_1,\ldots,x_n)$, be
a $K$-analytic function on an open and compact subset $X$ of
$K^n$. Let $g(y)$, $y=(y_1,\ldots,y_m)$, be a $K$-analytic
function on an open and compact subset $Y$ of $K^m$. Then
$f(x)+g(y)$ is a $K$-analytic function on the open and compact
subset $X \times Y$ of $K^{n+m}$. Put $A(s,\rho):=q^{s+1}-1$ if
$\rho$ is the trivial character of $R^{\times}$ and $A(s,\rho):=1$
if $\rho$ is another character of $R^{\times}$.

Fix a character $\chi$ of $R^{\times}$. Suppose that the only
critical value of $f$ and $g$ is zero. Then the poles of
$A(s,\chi) Z_{f+g,\chi}(s)$ are of the form $s_1+s_2$ with $s_1$ a
pole of $A(s,\chi')Z_{f,\chi'}(s)$ and $s_2$ a pole of
$A(s,\chi'')Z_{g,\chi''}(s)$ for some characters $\chi'$ and
$\chi''$ of $R^{\times}$ satisfying $\chi'\chi''=\chi$ (see
\cite{Igusalectures} or \cite[(5.1)]{Denefreport}).

\vspace{0,5cm}

\noindent \textbf{(4.2.1.2)} \textbf{Proposition.} \textsl{The set
$\{ s_0 \mid \exists f \in K\!\!<\!\!<\!\!x,y,z\!\!>\!\!> \, : \,
\mbox{mult}(f)=2$ and $Z_f(s)$ has a pole with real part $s_0 \}
\,  \cap \, ]-\infty,-1[$ is equal to
\[
\left\{ \left. -1-\frac{1}{i} \, \right| \, i \in \mathbb{Z}_{>1}
\right\}.
\] }


\noindent \textsl{Proof.} Let $f$ be an element of
$K\!\!<\!\!<\!\!x,y,z\!\!>\!\!>$ with multiplicity $2$. Up to an
affine coordinate transformation, the part of degree $2$ of $f$ is
equal to $ax^2+by^2+cz^2$, with $a,b,c \in K$ and $a \not= 0$.
Using (4.1.5), we may suppose that $f$ is of the form $x^2+g(y,z)$
with $g(y,z) \in K\!\!<\!\!<\!\!y,z\!\!>\!\!>$. The statement in
(4.2.1.1) and the result for curves imply that the real part of a
pole of $Z_f(s)$ is of the form $-1-1/i$, $i \in \mathbb{Z}_{>1}$,
if it is less than $-1$.

Now we prove the other inclusion. Using the $p$-adic stationary
phase formula \cite[Theorem 10.2.1]{Igusabook}, we obtain that
Igusa's $p$-adic zeta function of $xy+z^i$, $i \geq 2$, is equal
to
\[
\left( \frac{q-1}{q} \right) \left(
\frac{1-q^{-s-3}+(q-1)(q^{-2s-4}+q^{-3s-5}+ \cdots +
q^{-(i-1)s-(i+1)})}{(1-q^{-s-1})(1-q^{-is-(i+1)})} \right).
\]
The real poles of this zeta function are $-1$ and $-1-1/i$.
$\qquad \Box$

\subsubsection{Multiplicity larger than 2}
\noindent \textbf{(4.2.2.1)} Let $f$ be an element of
$K\!\!<\!\!<\!\!x,y,z\!\!>\!\!>$. Fix a (small enough)
neighbourhood $X$ of $0 \in K^3$ on which $f$ is convergent and an
embedded resolution $g:Y \rightarrow X$ of $f$ which is a
$K$-bianalytic map at the points where $f$ has normal crossings
and which is a composition of blowing-ups $g_{ij}:X_i \rightarrow
X_j$ with centre a $K$-analytic closed submanifold $D_{j}$ and
with exceptional surface $E_i$ satisfying:
\\ \indent (a) the codimension of $D_j$ in $X_j$ is at least 2; \\
\indent (b) $D_j$ is a subset of the zero locus of the strict
transform of $f$ on each chart (the strict transform of $f$ is not
defined globally);
\\ \indent (c) the union of the exceptional varieties in $X_j$ has
only normal crossings with $D_j$, i.e., for all $P \in D_j$, there
are three surface germs through $P$ which are in normal crossings
such that each exceptional surface germ through $P$ is one of them
and such that the germ of $D_j$ at $P$ is the intersection of some
of them;
\\ \indent (d) the image of $D_j$ in $X \subset K^3$ contains the
origin of $K^3$; and
\\ \indent (e) $D_j$ contains a point in which the pullback of
$f$ has not normal crossings. \\ Remark that such a resolution
always exists by Hironaka's theorem \cite{Hironaka}.

\vspace{0,5cm}

\noindent \textbf{(4.2.2.2)} The following table gives the
numerical data of $E_i$. In the columns, the dimension of $D_j$ is
kept fixed. In the rows, the number of exceptional surfaces
through $D_j$ is kept fixed. So $E_k$, $E_l$ and $E_m$ represent
exceptional surfaces that contain $D_j$. The multiplicity of the
strict transform of $f$ at $D_j$ is denoted by $\mu_{D_j}$.
\begin{center}
\begin{tabular}{ccc}
 \hline \vspace{-3mm} \\ & $D_j$ is a point $P$ & $D_j$ is a curve $L$ \\
 \hline \vspace{-3mm} \\ / & $(\mu_P,3)$ & $(\mu_L,2)$ \\
 $E_k$ & $(N_k+\mu_P,\nu_k+2)$ & $(N_k+\mu_L,\nu_k+1)$ \\
 $E_k$ and $E_l$ & $(N_k+N_l+\mu_P,\nu_k+\nu_l+1)$ &
 $(N_k+N_l+\mu_L,\nu_k+\nu_l)$ \\
 $E_k$, $E_l$ and $E_m$ & $(N_k+N_l+N_m+\mu_P,\nu_k+\nu_l+\nu_m)$
 & / \\ \hline
\end{tabular}
\end{center}

\vspace{0.5cm}

\noindent \textbf{(4.2.2.3) Lemma.} \textsl{Suppose that} mult$(f)
\geq 3$. \textsl{If there is no exceptional surface through $D_j$,
then $-\nu_i/N_i \geq -1$.}

\vspace{0,2cm} \noindent \emph{Proof.} The analogous statement for
the local topological zeta function is treated in
\cite[(3.3.3)]{SegersVeys}. The proof of the lemma is a trivial
adaptation of the proof there. $\qquad \Box$

\vspace{0,5cm}

\noindent \textbf{(4.2.2.4)} Suppose that $D_j$ is contained in at
least one exceptional surface and that the real parts of the
candidate poles associated to the exceptional surfaces that pass
through $D_j$ are larger than or equal to $-1$. Then the table in
(4.2.2.2) implies that also $-\nu_i/N_i \geq -1$, unless $D_j$ is
a regular point $P$ of the strict transform of $f$ around $P$
through which only one exceptional surface $E_0$ passes and
$-\nu_0/N_0=-1$. Suppose that we are in this situation. Let $Z_0$
be a (small enough) neighbourhood of $P$ such that, if we restrict
the blowing-ups $g_{ij}$ to the inverse image of $Z_0$, we get an
embedded resolution $h=h_1 \circ \cdots \circ h_s$ of the pullback
of $f$ which is a composition of blowing-ups $h_i:Z_i \rightarrow
Z_{i-1}$, $i \in \{1,\ldots,s\}$, with centre $D_{i-1}':=D_{i-1}
\cap Z_{i-1}$ and exceptional surface $E_i':=E_i \cap Z_i$ for
which $P$ is in the image of $D_{i-1}'$ under $h_1 \circ \cdots
\circ h_{i-1}$.

\vspace{0,2cm} \noindent Remark that it can happen that $g_{ij}$
is a $K$-bianalytic map on the inverse image of $Z_0$. Because we
did not specify the indices in (4.2.2.1), we were able to get a
nice notation here. From now on, we study the resolution $h:Z_s
\rightarrow Z_0$ of the pullback of $f$.

\vspace{0,5cm}

\noindent \textbf{Lemma.} \textsl{(a) If $D_i=D_i'$, then $D_i$ is
a subset of $E_0':=E_0 \cap Z_0$.}

\textsl{(b) Suppose that} mult$(f) \geq 3$. \textsl{Then we have}
$\nu_i \leq N_i+1$ \textsl{for every exceptional surface} $E_i$,
$i \in \{1,\ldots,s\}$. \textsl{Moreover,} $\nu_i=N_i+1$
\textsl{if and only if} $D_{i-1}$ \textsl{is a point and the
numerical data of every exceptional surface} $E_j$
\textsl{different from} $E_0$ \textsl{and through} $D_{i-1}$
\textsl{satisfy} $\nu_j=N_j+1$.

\textsl{(c) If} mult$(f) \geq 3$ \textsl{and if the numerical data
of} $E_i$ \textsl{satisfy} $\nu_i=N_i+1$, \textsl{then}
$-\nu_i/N_i \not= -\nu_j/N_j$ \textsl{for every exceptional
surface} $E_j$ \textsl{that intersects} $E_i$ \textsl{at some
stage of the resolution process.}

\vspace{0,2cm}

\noindent \textsl{Proof.} See \cite[(3.3.5),(3.3.6) and
(3.3.7)]{SegersVeys}. $\qquad \Box$

\vspace{0,5cm}

\noindent \textbf{Proposition.} \textsl{If} mult$(f) \geq 3$,
\textsl{then $Z_f(s)$ has no pole with real part less than $-1$.}

\vspace{0,2cm} \noindent \emph{Proof.} The proof is analogous to
the one in \cite[(3.3.8)]{SegersVeys}, except that we have to use
the proposition in 3.2.  $\qquad \Box$

\vspace{1cm}

\noindent {\Large \textbf{Appendix. Poles and divisibility of
the}} \textbf{$\emph{{\Large M}}_\emph{i}$}

\vspace{0,5cm}

Suppose that $f$ is a $K$-analytic function on $R^n$ defined by a
power series over $R$ which is convergent on the whole of $R^n$.
Let $l$ be the smallest real part of a pole of $Z_f(s)$ and let
$M_i$ be the number of solutions of $f(x) \equiv 0 \mbox{ mod }
P^i$ in $(R/P^i)^n$.

\vspace{0,2cm}

\noindent \textbf{Proposition.} \textsl{There exists an integer
$a$ which is independent of $i$ such that $M_i$ is an integer
multiple of $q^{\ulcorner(n+l)i-a\urcorner}$ for all $i \in
\mathbb{Z}_{\geq 0}$.}

\vspace{0,2cm}

\noindent \textsl{Remark.} (i) The number
$\ulcorner(n+l)i-a\urcorner$ is the smallest integer larger than
or equal to $(n+l)i-a$, which rises ($n+l > 0$) linearly as a
function of $i$ with a slope depending on $l$. \\ (ii) The
statement is trivial if $(n+l)i-a \leq 0$ because the $M_i$ are
integers. If $(n+l)i-a>0$, which is the case for $i$ large enough,
it claims that $M_i$ is divisible by
$q^{\ulcorner(n+l)i-a\urcorner}$.

\vspace{0,2cm}

\noindent \textsl{Proof.} Put $t=q^{-s}$. It follows from (2.2)
that we can write
\[
Z_f(t) = \frac{A(t)}{\prod_{j \in T} \left( 1- q^{-\nu_j}t^{N_j}
\right) },
\]
where $A(t)$ is a polynomial with coefficients in the set $S:=\{
z/q^i \mid z \in \mathbb{Z}, i \in \mathbb{Z}_{\geq 0} \}$. By
using the division algorithm for polynomials we can write
\[
Z_f(t) = \frac{B(t)}{\prod_{j \in K} \left( 1- q^{-\nu_j}t^{N_j}
\right) },
\]
where $B(t)$ is a polynomial with coefficients in $S$ and where
$K:=\{ j \in T \mid -\nu_j/N_j \geq l \}$.

The Poincar\'e series $P(t)$ of $f$ is defined by
\[
P(t) = \sum_{i=0}^{\infty} M_i \frac{t^i}{q^{ni}}.
\]
and can be obtained from $Z_f(t)$ by the relation
\[
P(t) = \frac{1-tZ_f(t)}{1-t}.
\]
It easily follows from the defining integral of Igusa's $p$-adic
zeta function that $Z_f(t=1)=1$. Consequently, $1-tZ_f(t)$ is
divisible by $1-t$ and $P(t)$ can be written as
\[
P(t) = \frac{C(t)}{\prod_{j \in K} \left( 1- q^{-\nu_j} t^{N_j}
\right) },
\]
where $C(t)$ is a polynomial with coefficients in $S$.

We will say that a formal power series in $t$ has the divisibility
property if the coefficient of $t^i/q^{ni}$ is an integer multiple
of $q^{\ulcorner(n+l)i\urcorner}$ for every $i$.

For $j \in K$, the series
\[
\frac{1}{1-q^{-\nu_j}t^{N_j}}= \sum_{i=0}^{\infty} q^{-i\nu_j}
t^{iN_j} = \sum_{i=0}^{\infty} q^{i(nN_j-\nu_j)}
\frac{t^{iN_j}}{q^{niN_j}}
\]
has the divisibility property because $nN_j-\nu_j$ is an integer
larger than or equal to $N_j(n+l)$. Let $a$ be an integer such
that the polynomial $D(t):=q^a C(t)$ has the divisibility
property. Remark that $C(t)=q^{-a} D(t)$.

One can easily check that the product of a finite number of power
series with the divisibility property also has the divisibility
property. This implies that $P(t)$ is a power series with the
divisibility property, multiplied by $q^{-a}$. Hence $M_i$ is an
integer multiple of $q^{\ulcorner(n+l)i\urcorner -a} =
q^{\ulcorner(n+l)i-a\urcorner}$ for all $i$. $\qquad \Box$

\footnotesize{

\noindent \textsc{K.U.Leuven, Departement Wiskunde,
Celestijnenlaan 200B, B-3001 Leuven, Belgium,} \\ \textsl{E-mail:}
dirk.segers@wis.kuleuven.be} \\ \textsl{URL:}
http://wis.kuleuven.be/algebra/segers/segers.htm

\end{document}